\documentclass[12pt,a4paper]{article}

\usepackage{epsfig}
\usepackage{amsmath}
\usepackage{graphicx}
\usepackage{geometry}         
\usepackage{latexsym}
\usepackage[latin1]{inputenc} 
\usepackage{amsmath}

\newtheorem{theorem}{Theorem}
\newtheorem{lemma}{Lemma}

\title{Solvability of the boundary value problem associated with the wave diffraction by a layer filled with a Kerr-type nonlinear medium}

\author{Yury Shestopalov~$^1$ and  Vasil Yatsyk~$^2$}

\date{$^1$ Karlstad University, Karlstad, Sweden \\~~E-mail:
youri.shestopalov@kau.se\\[10pt]
$^2$ The Usikov Institute of Radio Physics and Electronics,
Kharkov, Ukraine\\~~E-mail: yatsyk@vk.kharkov.ua}

\begin{document}

\maketitle

\tableofcontents


\begin{abstract}
\noindent
 The diffraction of a plane wave by a transversely
inhomogeneous isotropic nonmagnetic linearly polarized dielectric
layer filled with a Kerr-type nonlinear medium is considered. The
analytical and numerical solution techniques are developed. The
diffraction problem is reduced to a singular boundary value
problem for a semilinear second-order ordinary differential
equation with a cubic nonlinearity and then to a cubic-nonlinear
integral equation of the second kind and to a system of nonlinear
operator equations of the second kind solved using iterations.
Sufficient conditions of the unique solvability are obtained using
the contraction principle.
\end{abstract}




\section{Introduction } \label{subsec:mylabel1}

 Scattering and propagation of electromagnetic
waves in layered structures filled with nonlinear media have been
a subject of intense studies since the 1970s. A goal of this work
is to develop solution techniques for singular boundary value
problems (BVPs) for the Maxwell equations arising in mathematical
models of the wave diffraction in nonlinear media elaborated in
\cite{[1]}--\cite{[7]},
\cite{[3]}, \cite{[4]}, \cite{[5]}
 that can be reduced to one-dimensional settings for the
Helmholtz equation on the line \cite{[1]}, \cite{[2]}. The BVPs
are formulated on infinite and semi-infinite intervals and with
transmission-type conditions and conditions at infinity that
contain the spectral parameter \cite{[3]}, \cite{[4]}, \cite{[5]}.
When the wave propagation in a cylindrical dielectric waveguide
filled with a nonlinear medium is considered \cite{[4]},
\cite{[5]}, the coefficient in the equation multiplying the
nonlinear term differs from zero inside a finite interval $\left(
{0,\,a} \right)$ and the conditions are stated at the point $a$
(continuity), at the origin (e.g. boundedness), and at infinity
(rate of decay). The corresponding singular semilinear BVPs are
formulated for the differential operators $L\left( \lambda
\right)\,u + \alpha \,B\left( {u{\kern 1pt} ;\,\lambda } \right) =
0$, where $L{\kern 1pt} u$ is a linear differential operator and
$B\left( {u{\kern 1pt} ;\,\lambda } \right)$ is a nonlinear
operator. An example with $B\left( {u{\kern 1pt} } \right) = u^3$
associated with the study of the wave propagation in Kerr-type
nonlinear fibers is considered in great detail in  \cite{[4]},
\cite{[5]}. The method of solution employs reduction to nonlinear
integral equations (IEs)  \cite{[4]}, \cite{[5]}, \cite{[6]},
\cite{[7]}, \cite{[8]} constructed using Green's function of the
linear differential operator $L{\kern 1pt} u$; the eigenvalue
problems are then replaced by the determination of characteristic
numbers of integral operator-valued functions (OVFs) that are
nonlinear both with respect to the solution and the spectral
parameter. The latter problems are reduced to the functional
dispersion equations, and their roots give the sought-for
eigenvalues. The existence and distribution of roots on the
complex plane are verified. The linearization is considered in
\cite{[6]}.

The reflection and transmission of electromagnetic waves at a
nonlinear homogeneous, isotropic, non-magnetic dielectric layer
situated between two linear homogeneous, semi-infinite media is of
particular interest in linear optics \cite{[9]}. In nonlinear
optics, the Kerr-like nonlinear dielectric film has been the focus
of a number of studies \cite{[10]}, \cite{[11]}, \cite{[12]},
\cite{[13]}, \cite{[14]}. In \cite{[10]}, \cite{[11]},
\cite{[12]}, \cite{[13]}, the solutions of the nonlinear Helmholtz
equations have been given in terms of various Jacobian elliptic
functions. The explicit form of these functions depends on the
associated parameter regimes. As shown in \cite{[14]},
 no classification of the
solutions with respect to different parameter regimes is
necessary, since the general solution can be presented in terms of
Weierstrass' elliptic functions containing the complete parameter
dependence. In \cite{[1]}, a simplified version of this result is
given, generalizing the approach applied in linear optics. Namely,
a general analytical solution of the Helmholtz equation is
obtained describing the scattering of a plane, monochromatic,
TE-polarized wave by a transversely homogeneous dielectric layer
(with a constant permittivity) exhibiting a local Kerr-like
nonlinearity. The layer is situated between two semi-infinite
non-absorbing, non-magnetic, isotropic, and homogeneous media. The
results derived contain the conditions for unbounded field
intensities expressed in terms of the imaginary half-period of
Weierstrass' elliptic function. The reflectivity $R$ is calculated
as a function of the layer thickness and the transmitted
intensity. The critical values of $R$ are determined.

In \cite{[7]} the approach set forth in \cite{[1]}, \cite{[2]} is
applied to the analysis of the problems of the wave diffraction by
layers filled both with linear and nonlinear dielectric media
having constant and variable permittivities. The plane wave
diffraction problem is reduced in \cite{[7]} to a nonlinear
Volterra IE and its solution is obtained as a limit of a certain
function sequence. The sufficient conditions for the IE unique
solvability are obtained by estimating the norms of the associated
nonlinear operator.

In this paper the approaches developed in \cite{[4]}, \cite{[5]}
and \cite{[3]}, \cite{[15]} \cite{[16]} are applied to the
solution of singular semilinear BVPs arising in a mathematical
model of the wave diffraction from a transversely inhomogeneous
dielectric layer having a variable permittivity. The approach
employs Fredholm-type IEs with complex-valued kernels derived on
the basis of the method proposed in \cite{[3]} and differs thus
from the technique \cite{[1]}, \cite{[7]}, \cite{[8]} based on the
reduction to nonlinear Volterra IEs. On the other hand, the
sufficient solvability conditions presented in this study are
different from those reported in \cite{[7]}; in fact, these
conditions are obtained explicitly in terms of the problem
parameters. Next, in this paper we apply the solution technique
based on the analysis of cubic-nonlinear IEs to prove the unique
solvability of the diffraction problem for a lossy weakly
nonlinear layer with a complex-valued permittivity function. We
note in this respect \cite{[8]} where this problem is solved for a
layer filled by linear and nonlinear lossy media using a general
approach which enables one to evaluate the solutions in terms of
uniformly convergent sequences of iterations of the Volterra IEs.
\section{Maxwell equations and wave propagation in nonlinear
media}

\subsection{General assumptions leading to the problem statement}
%
Nonlinear processes in electrodynamics and optics are
described by the Maxwell equations
\begin{equation}
\label{eq1}
\begin{array}{l}
 \nabla \times \vec {E}\left( {\vec {r},t} \right) = -
\frac{1}{c}\frac{\partial \vec {B}\left( {\vec {r},t} \right)}{\partial
\,t}\,,\quad \nabla \times \vec {H}\left( {\vec {r},t} \right) =
\frac{1}{c}\frac{\partial \vec {D}\left( {\vec {r},t} \right)}{\partial
\,t}\,, \\
 \nabla \cdot \vec {D}\left( {\vec {r},t} \right) = 0\,,\quad \quad \quad
\quad \quad \;\,\nabla \cdot \vec {B}\left( {\vec {r},t} \right) = 0, \\
 \end{array}
\end{equation}
Here $\vec {E}\left( {\vec {r},t} \right)$, $\vec {H}\left( {\vec
{r},t} \right) \quad \vec {D}\left( {\vec {r},t} \right)$, and
$\vec {B}\left( {\vec {r},t} \right)$ are the vectors of,
respectively, electric and magnetic field intensities, electric
displacement, and magnetic induction. This system is complemented
by material equations

\begin{equation}
\label{eq2}
\begin{array}{l}
 \vec {D}\left( {\vec {r},t} \right) = \vec {E}\left( {\vec {r},t} \right) +
4\pi \vec {P}\left( {\vec {r},t} \right), \\
 \vec {B}\left( {\vec {r},t} \right) = \vec {H}\left( {\vec {r},t} \right) +
4\pi \vec {M}\left( {\vec {r},t} \right), \\
 \end{array}
\end{equation}

\noindent where $\vec {P}\left( {\vec {r},t} \right)$ and $\vec
{M}\left( {\vec {r},t} \right)$ are the vectors of, respectively,
polarization and magnetic moment.

The polarization vector $\vec {P}\left( {\vec {r},t} \right) =
\hat {F}\left[ {\vec {E}\left( {\vec {r},t} \right)} \right]$,
where $\hat {F}$ denotes a certain nonlinear operator, is
generally nonlinear (with respect to the intensity) and nonlocal
both in time and space In this work, we will limit the analysis,
following \cite{[17]}, to nonlinear media having spatially
nonlocal response function. In this case the polarization vector
can be expanded \cite{[18]} in terms of the electric field
components

\begin{equation}
\label{eq3}
P_i \left( {\vec {r},t} \right) \equiv \chi _{ij}^{\left( 1 \right)} E_j +
\chi _{ijk}^{\left( 2 \right)} E_j E_k + \chi _{ijkl}^{\left( 3 \right)} E_j
E_k E_l + \ldots
\end{equation}

Here $P_i $ and $E_i $ are the components of the polarization and
electric vectors, respectively and coefficients $\chi $ are lower
terms of the expansion for nonlinear susceptibility. Thus, we
ignore spatial dispersion \cite{[19]}. Note however that this
assumption does not limit our possibilities to consider effects in
media with interfaces where linear and nonlinear susceptibility
tensors may depend on spatial variables \cite{[17]}.

Below, we assume that the medium is nonmagnetic, $\vec {M}\left( {\vec
{r},t} \right) \equiv 0$. Resolving equations (\ref{eq1}) and (\ref{eq2}) with respect to
$\vec {H}\left( {\vec {r},t} \right)$ we reduce them to one vector equation

\begin{equation}
\label{eq4}
\nabla ^2\vec {E}\left( {\vec {r},t} \right) - \nabla \left[ {\nabla \cdot
\vec {E}\left( {\vec {r},t} \right)} \right] - \frac{1}{c^2}\frac{\partial
^2}{\partial \,t^2}\vec {D}^{\left( L \right)}\left( {\vec {r},t} \right) -
\frac{4\pi }{c^2}\frac{\partial ^2}{\partial \,t^2}\vec {P}^{\left( {NL}
\right)}\left( {\vec {r},t} \right) = 0,
\end{equation}

\noindent
where $\vec {D}^{\left( L \right)} = \vec {E} + 4\pi \vec {P}^{\left( L
\right)} = \hat {\varepsilon }\vec {E}$, $\vec {P}^{\left( L \right)} = \hat
{\chi }^{\left( 1 \right)}\vec {E}$, and $\hat {\varepsilon } = \varepsilon
^{\left( L \right)}$ are the linear terms of the electric displacement and
polarization vectors and permittivity tensor, respectively (here
$D_i^{\left( L \right)} = \varepsilon _{ij}^{\left( L \right)} E_j $, $P_i
^{\left( L \right)} \equiv \chi _{ij}^{\left( 1 \right)} E_j $, and
$\varepsilon _{ij}^{\left( L \right)} = 1 + 4\pi \chi _{ij}^{\left( 1
\right)} )$; $\vec {P}^{\left( {NL} \right)}$ is the nonlinear part of the
polarization vector (according to (\ref{eq3}), $P_i ^{\left( {NL} \right)} \equiv
\chi _{ijk}^{\left( 2 \right)} E_j E_k + \chi _{ijkl}^{\left( 3 \right)} E_j
E_k E_l + \ldots )$; and$\chi _{ij}^{\left( 1 \right)} $, $\chi
_{ijk}^{\left( 2 \right)} $, $\chi _{ijkl}^{\left( 3 \right)} $, are the
respective components of the medium susceptibility tensors $\hat {\chi
}^{\left( 1 \right)}$, $\hat {\chi }^{\left( 2 \right)}$, $\hat {\chi
}^{\left( 3 \right)}$.

Equation (\ref{eq4}) is of general character and is used, together
with material equations (\ref{eq2}), in electrodynamics and
optics; in every particular case, specific assumptions are made
that enable one to simplify its form. Note, for example, that in
the majority of important problems the longitudinal field
components (along the $z$-axis) are negligible \cite{[17]}. The
second term in (\ref{eq4}) $\nabla \left[ {\nabla \cdot \vec
{E}\left( {\vec {r},t} \right)} \right]$ (where the inner product
can be written, using the condition$\nabla \cdot \vec {D} = 0$, in
the form$\nabla \cdot \vec {E} = - \,\left[ {{\vec {E} \cdot
\left( {\nabla \hat {\varepsilon }} \right)} \mathord{\left/
{\vphantom {{\vec {E} \cdot \left( {\nabla \hat {\varepsilon }}
\right)} {\hat {\varepsilon }\,}}} \right.
\kern-\nulldelimiterspace} {\hat {\varepsilon }\,}} \right]\,)$
contains both longitudinal and transverse field components and may
be ignored in a number of cases.

Consider the case of stationary electromagnetic field

\begin{equation}
\nonumber \vec {E}\left( {\vec {r},t} \right) = Re\left[ {\exp
\left( { - i\omega \,t} \right)  \vec {E}\left( \vec {r} \right)}
\right] \equiv \frac{1}{2}\left[ {\exp \left( { - i\omega \,t}
\right) \vec {E}\left( \vec {r} \right) + \exp \left( {i\omega
\,t} \right)  \vec {E}^\ast \left( \vec {r} \right)} \right]\, ,
\end{equation}

\begin{equation}
\nonumber \vec {H}\left( {\vec {r},t} \right) = Re\left[ {\exp
\left( { - i\omega \,t} \right)  \vec {H}\left( \vec {r} \right)}
\right] \equiv \frac{1}{2}\left[ {\exp \left( { - i\omega \,t}
\right)  \vec {H}\left( \vec {r} \right) + \exp \left( {i\omega
\,t} \right)  \vec {H}^\ast \left( \vec {r} \right)} \right]\, ,
\end{equation}

\noindent where $\vec {E}\left( \vec {r} \right)$ and $\vec
{H}\left( \vec {r} \right)$ are the complex amplitudes of the
electric, $\vec {E}\left( {\vec {r},t} \right)$, and magnetic,
$\vec {H}\left( {\vec {r},t} \right)$, field intensity vectors,
$Re$ is the real part of the complex vector-function, and $\,\ast
\,$denotes complex conjugation. Assuming that the medium is weakly
nonlinear (when the so-called weakly-waveguide approximation
holds), i.e.
\begin{equation}
\label{eq5}
\left| {\varepsilon _{ij}^{\left( {NL} \right)} } \right| < < \left|
{\varepsilon _{ij}^{\left( L \right)} } \right|,
\end{equation}
where $\varepsilon _{ij}^{\left( {NL} \right)} = 4\pi \chi
_{ijkl}^{\left( 3 \right)} E_k E_l^\ast $ is governed by nonlinear
terms in (\ref{eq4}) (see \cite{[17]}, \cite{[20]}), one can
generalize weakly-waveguide approximation \cite{[21]}  and take
into account the effect of nonlinear self-canalization
\cite{[22]}. In this case one can ignore the second term in
(\ref{eq4}), which is equivalent to ignoring the longitudinal
field components, and vectors $\vec {E}$ and $\vec {P}$ will have
only transverse components \cite{[17]}.

%

Consider the diffraction of a stationary electromagnetic wave
$\left[ {\sim \exp \left( { - i\omega \,t} \right)} \right]$ by a
weakly nonlinear object. Perform a transition to the frequency
domain using the direct and inverse Fourier transforms
\begin{equation}
\nonumber
 \left[ {{\begin{array}{*{20}c}
 {\dot {\vec {E}}\left( {\vec {r},\tilde {\omega }} \right)} \hfill \\
 {\dot {\vec {D}}^{\left( L \right)}\left( {\vec {r},\tilde {\omega }}
\right)} \hfill \\
 {\dot {\vec {P}}^{\left( {NL} \right)}\left( {\vec {r},\tilde {\omega }}
\right)} \hfill \\
\end{array} }} \right] = \int\limits_{ - \infty }^\infty {\left[
{{\begin{array}{*{20}c}
 {\vec {E}\left( {\vec {r},t} \right)} \hfill \\
 {\vec {D}^{\left( L \right)}\left( {\vec {r},t} \right)} \hfill \\
 {\vec {P}^{\left( {NL} \right)}\left( {\vec {r},t} \right)} \hfill \\
\end{array} }} \right]\,e^{i{\kern 1pt} \tilde {\omega }{\kern 1pt} t}dt}
, \, \left[ {{\begin{array}{*{20}c}
 {\vec {E}\left( {\vec {r},t} \right)} \hfill \\
 {\vec {D}^{\left( L \right)}\left( {\vec {r},t} \right)} \hfill \\
 {\vec {P}^{\left( {NL} \right)}\left( {\vec {r},t} \right)} \hfill \\
\end{array} }} \right] = \frac{1}{2\pi }\int\limits_{ - \infty }^\infty
{\left[ {{\begin{array}{*{20}c}
 {\dot {\vec {E}}\left( {\vec {r},\tilde {\omega }} \right)} \hfill \\
 {\dot {\vec {D}}^{\left( L \right)}\left( {\vec {r},\tilde {\omega }}
\right)} \hfill \\
 {\dot {\vec {P}}^{\left( {NL} \right)}\left( {\vec {r},\tilde {\omega }}
\right)} \hfill \\
\end{array} }} \right]\,e^{ - i{\kern 1pt} \tilde {\omega }{\kern 1pt}
t}d\tilde {\omega }} .
\end{equation}

Applying formally the Fourier transform to equation (\ref{eq4}) we
obtain the following representation in the frequency domain
\begin{equation}
\label{eq6}
\nabla ^2\dot {\vec {E}}\left( {\vec {r},\tilde {\omega }} \right) - \nabla
\left[ {\nabla \cdot \dot {\vec {E}}\left( {\vec {r},\tilde {\omega }}
\right)} \right] + \frac{\omega ^2}{c^2}\dot {\vec {D}}^{\left( L
\right)}\left( {\vec {r},\tilde {\omega }} \right) + \frac{4\pi \omega
^2}{c^2}\dot {\vec {P}}^{\left( {NL} \right)}\left( {\vec {r},\tilde {\omega
}} \right) = 0.
\end{equation}

A stationary $\left[ {\sim \exp \left( { - i\omega \,t} \right)}
\right]$ electromagnetic wave propagating in a weakly nonlinear
dielectric structure gives rise to a field containing all
frequency harmonics, see \cite{[19]}, \cite{[20]}.
Therefore, the quantities describing the electromagnetic field in
the time domain subject to equation (\ref{eq4}) can be represented
as Fourier series
\begin{equation}
\label{eq7}
\begin{array}{l}
 \vec {E}\left( {\vec {r},t} \right) = \frac{1}{2}\sum\limits_{n = - \infty
}^\infty {\vec {E}\left( {\vec {r},n\omega } \right)\exp \left( {
- i n \omega t} \right)} ,\quad \vec {D}^{\left( L \right)}\left(
{\vec {r},t} \right) = \frac{1}{2}\sum\limits_{n = - \infty
}^\infty {\vec {D}^{\left( L \right)}\left( {\vec {r},n\omega }
\right)\exp \left( { - i n \omega t} \right)\,} ,\; \\
 \vec {P}^{\left( {NL} \right)}\left( {\vec {r},t} \right) =
\frac{1}{2}\sum\limits_{n = - \infty }^\infty {\vec {P}^{\left( {NL}
\right)}\left( {\vec {r},n\omega } \right)\exp \left( { - i n \omega t} \right)\,} .
 \end{array}
\end{equation}

Applying to (\ref{eq7}) the Fourier transform we obtain
\begin{equation}
\label{eq8}
\begin{array}{l}
 \left[ {{\begin{array}{*{20}c}
 {\dot {\vec {E}}\left( {\vec {r},\tilde {\omega }} \right)} \hfill \\
 {\dot {\vec {D}}^{\left( L \right)}\left( {\vec {r},\tilde {\omega }}
\right)} \hfill \\
 {\dot {\vec {P}}^{\left( {NL} \right)}\left( {\vec {r},\tilde {\omega }}
\right)} \hfill \\
\end{array} }} \right] = \int\limits_{ - \infty }^\infty {\left[
{{\begin{array}{*{20}c}
 {\vec {E}\left( {\vec {r},t} \right)} \hfill \\
 {\vec {D}^{\left( L \right)}\left( {\vec {r},t} \right)} \hfill \\
 {\vec {P}^{\left( {NL} \right)}\left( {\vec {r},t} \right)} \hfill \\
\end{array} }} \right]\,e^{i\tilde{\omega } t}\,dt}
= \\
 = \frac{1}{2}\int\limits_{ - \infty }^\infty {\;\sum\limits_{n = - \infty
}^\infty {\left[ {{\begin{array}{*{20}c}
 {\vec {E}\left( {\vec {r},n\omega } \right)} \hfill \\
 {\vec {D}^{\left( L \right)}\left( {\vec {r},n\omega } \right)} \hfill \\
 {\vec {P}^{\left( {NL} \right)}\left( {\vec {r},n\omega } \right)} \hfill
\\
\end{array} }} \right]\,e^{ - i n {\omega } t} e^{i\tilde{\omega } t}\,dt} = \frac{\sqrt {2\pi }
}{2} {\kern 1pt} \left[ {{\begin{array}{*{20}c}
 {\vec {E}\left( {\vec {r},n\omega } \right)} \hfill \\
 {\vec {D}^{\left( L \right)}\left( {\vec {r},n\omega } \right)} \hfill \\
 {\vec {P}^{\left( {NL} \right)}\left( {\vec {r},n\omega } \right)} \hfill
\\
\end{array} }} \right]{\kern 1pt} \delta \left. {\left( 0 \right)}
\right|_{\tilde {\omega } = n\omega }} \,, \\
 \end{array}
\end{equation}

\noindent
where $\delta \left( s \right) = \frac{1}{\sqrt {2\pi } }\int\limits_{ -
\infty }^\infty {\exp \left( {i{\kern 1pt} s{\kern 1pt} t} \right)dt} $ is
the Dirac delta-function.

Substituting (\ref{eq8}) into (\ref{eq6}), we obtain an infinite equation system with
respect to the sought-for Fourier amplitudes of the electromagnetic of the
weakly nonlinear structure in the frequency domain,
\begin{eqnarray}
\label{eq9} \nabla ^2\vec {E}\left( {\vec {r},n\omega } \right) -
\nabla \left[ {\nabla \cdot \vec {E}\left( {\vec {r},n\omega }
\right)} \right] + \frac{\omega ^2}{c^2}\vec {D}^{\left( L
\right)}\left( {\vec {r},n\omega } \right) + \frac{4\pi \omega
^2}{c^2}\vec {P}^{\left( {NL} \right)}\left( {\vec {r},n\omega }
\right) = 0,
\\ \nonumber n = 0,\,\pm 1,\,\pm 2,\,\ldots \, .
\end{eqnarray}

For linear electrodynamic objects the equations in the system
(\ref{eq9}) are independent. In a nonlinear structure, the
presence of functions $\vec {P}^{\left( {NL} \right)}\left( {\vec
{r},n\omega } \right)$ makes them coupled since every harmonic
depends on a series of $\vec {E}\left( {\vec {r},n\omega }
\right)$. Indeed consider a three-component electromagnetic field
$\vec {E} = \left( {E_x ,0,0} \right)$, $\vec {H} = \left( {0,H_y
,H_z } \right)$. The fact that the field $\vec {E} = \left( {E_x
,0,0} \right)$ has one component enables one to consider
(\ref{eq9}) as a system of scalar equations with respect to $E_x
$. Take lower terms in the expansion (\ref{eq3}) in the vicinity
of the zero value of the electric field intensity. Then the only
nonzero component of the polarization vector $\vec {P} = \left(
{P_x ,0,0} \right)$ is determined by the third-order
susceptibility tensor $\hat {\chi }^{\left( 3 \right)}$, which is
characteristic for the Kerr-type medium. In the time domain, this
component can be represented in the form (cf. (\ref{eq3}) and
(\ref{eq7})):

\begin{equation}
\label{eq10}
\begin{array}{l}
 P_x \left( {\vec {r},t} \right) = \frac{1}{2}\sum\limits_{s = - \infty
}^\infty {P_x \left( {\vec {r},s\omega } \right)\exp \left( { - i\omega
{\kern 1pt} s{\kern 1pt} t} \right)} \equiv \chi _{x{\kern 1pt} x{\kern 1pt}
x{\kern 1pt} x}^{\left( 3 \right)} E_x \left( {\vec {r},t} \right)E_x \left(
{\vec {r},t} \right)E_x \left( {\vec {r},t} \right) = \\
 = \frac{1}{8}\sum\limits_{\left\{ {{\begin{array}{*{20}c}
 {n,{\kern 1pt} m,{\kern 1pt} p,{\kern 1pt} s = - \infty } \hfill \\
 {n + m + p = s} \hfill \\
\end{array} }} \right.}^\infty {\chi _{x{\kern 1pt} x{\kern 1pt} x{\kern
1pt} x}^{\left( 3 \right)} \left( {s\omega ;\,n\omega ,\,m\omega
,\,p\omega } \right)E_x \left( {\vec {r},n\omega } \right)E_x
\left( {\vec {r},m\omega } \right)E_x \left( {\vec {r},p\omega }
\right) e^{- i\omega (n + m + p) t} }.
 \end{array}
\end{equation}

Applying to (\ref{eq10}) the Fourier transform with respect to time 
(\ref{eq8}) we obtain an expansion in the frequency domain

\begin{equation}
\label{eq11}
\begin{array}{l}
 P_x \left( {\vec {r},s\omega } \right) = \frac{1}{4}\sum\limits_{\left\{
{{\begin{array}{*{20}c}
 {n,{\kern 1pt} m,{\kern 1pt} p = - \infty } \hfill \\
 {n + m + p = s} \hfill \\
\end{array} }} \right.}^\infty {\chi _{x{\kern 1pt} x{\kern 1pt} x{\kern
1pt} x}^{\left( 3 \right)} \left( {s\omega \,;\,n\omega ,\,m\omega
,\,p\omega } \right)E_x \left( {\vec {r},n\omega } \right)E_x \left( {\vec
{r},m\omega } \right)E_x \left( {\vec {r},p\omega } \right)} = \\
 = \frac{1}{4}\sum\limits_{j = 0}^\infty {3\chi _{x{\kern 1pt} x{\kern 1pt}
x{\kern 1pt} x}^{\left( 3 \right)} \left( {s\omega \,;\,j\omega ,\, -
j\omega ,\,s\omega } \right)\left| {E_x \left( {\vec {r},j\omega } \right)}
\right|^2E_x \left( {\vec {r},s\omega } \right)} + \\
 + \frac{1}{4}\sum\limits_{\left\{ {{\begin{array}{*{20}c}
 {\begin{array}{l}
 n,{\kern 1pt} m,{\kern 1pt} p = - \infty \\
 n \ne - m,\,p = s \\
 m \ne - p,\,n = s \\
 n \ne - p,\,m = s \\
 \end{array}} \hfill \\
 {n + m + p = s} \hfill \\
\end{array} }} \right.}^\infty {\chi _{x{\kern 1pt} x{\kern 1pt} x{\kern
1pt} x}^{\left( 3 \right)} \left( {s\omega \,;\,n\omega ,\,m\omega
,\,p\omega } \right)E_x \left( {\vec {r},n\omega } \right)E_x \left( {\vec
{r},m\omega } \right)E_x \left( {\vec {r},p\omega } \right)} \,. \\
 \end{array}
\end{equation}

The addends in the first sum of (\ref{eq11}) are usually called
the phase self-modulation (PSM) terms \cite{[17]}. We obtained
them taking into account the property of the Fourier coefficients
$E_x \left( {\vec {r},j\omega } \right) = E_x^\ast \left( {\vec
{r},\, - j\omega } \right)$; factors 3 appear as a result of
permutations $\left\{ {j\omega ,\, - j\omega ,\,s\omega }
\right\}$ of three last parameters in the terms $\chi _{x{\kern
1pt} x{\kern 1pt} x{\kern 1pt} x}^{\left( 3 \right)} \left(
{s\omega \,;\,j\omega ,\, - j\omega ,\,s\omega } \right)$.

When particular nonlinear effects are considered, one can limit the analysis
to finitely many equations of system (\ref{eq9}), leaving in the formulas for the
polarization coefficients separate terms that characterize the physical
problem in question.

In this paper, we analyze electromagnetic fields scattered by a dielectric
layer filled by a Kerr-type (weakly) nonlinear medium. We limit the analysis
to such a level of intensities of the incident electromagnetic field
affecting the structure when harmonic oscillations at combined frequencies
may be neglected. In this case equations (\ref{eq9}) and (\ref{eq11}) have the form

\begin{equation}
\label{eq14}
\begin{array}{l}
 \nabla ^2E_x \left( {\vec {r},\omega } \right) - \nabla \left[ {\nabla
\cdot E_x \left( {\vec {r},\omega } \right)} \right] + \frac{\omega
^2}{c^2}D_x ^{\left( L \right)}\left( {\vec {r},\omega } \right) +
\frac{4\pi \omega ^2}{c^2}P_x ^{\left( {NL} \right)}\left( {\vec {r},\omega
} \right) = 0\,, \\
 P_x \left( {\vec {r},\omega } \right) = \frac{3}{4}\chi _{x{\kern 1pt}
x{\kern 1pt} x{\kern 1pt} x}^{\left( 3 \right)} \left( {\omega \,;\,\omega
,\, - \omega ,\,\omega } \right)\left| {E_x \left( {\vec {r},\omega }
\right)} \right|^2E_x \left( {\vec {r},\omega } \right)\,. \\
 \end{array}
\end{equation}

In the next sections, we will use equations (\ref{eq14}) to formulate a boundary value problem associated with the electromagnetic wave diffraction by a layer filled by a Kerr-type nonlinear medium and construct the methods of analytical solution to this problem.

\subsection{Statement of the problem of diffraction by a weakly nonlinear layer}
Denote by $\vec {E}\left( \vec {r} \right) \equiv \vec {E}\left(
{\vec {r},\omega } \right)$ and $\vec {H}\left( \vec {r} \right)
\equiv \vec {H}\left( {\vec {r},\omega } \right)$ the complex
amplitudes of the stationary electromagnetic field; the time
dependence is $\exp \left( { - i\omega \,t} \right)$. Consider the
problem of diffraction of a plane stationary electromagnetic wave
$\vec {E}\left( {\vec {r},t} \right) = \exp \left( { - i\omega
\,t} \right) \vec {E}\left( \vec {r} \right)$, $\vec {H}\left(
{\vec {r},t} \right) = \exp \left( { - i\omega \,t} \right) \vec
{H}\left( \vec {r} \right)$ by a nonmagnetic, $\vec {M} = 0$,
isotropic and linearly polarized $\vec {E}\left( \vec {r} \right)
= \left( {E_x \left( {y,z} \right),0,0} \right)$, $\vec {H}\left(
\vec {r} \right) = \left( {0,H_y \equiv \frac{1}{i\omega \mu _0
}\frac{\partial E_x }{\partial z},H_z \equiv - \frac{1}{i\omega
\mu _0 }\frac{\partial E_x }{\partial y}} \right)$
(E-polarization), transversely inhomogeneous, $\varepsilon
^{\left( L \right)}(z) = \varepsilon _{xx}^{\left( L \right)}
(z)$, dielectric layer with a weak Kerr-type nonlinearity
(\ref{eq5}) $P_x^{\left( {NL} \right)} = \left( {3 \mathord{\left/
{\vphantom {3 4}} \right. \kern-\nulldelimiterspace} 4}
\right)\chi _{xxxx}^{\left( 3 \right)} \left| {E_x } \right|^2E_x
$ (where $\vec {P}^{\left( {NL} \right)} = (P_x^{\left( {NL}
\right)} ,\;\;0,\;0)$;
this problem is stated in
\cite{[1]}, \cite{[7]}, \cite{[17]}, and \cite{[21]}. Using
(\ref{eq1}), (\ref{eq2}), and the results from \cite{[17]} we
obtain $\nabla \cdot \vec {E} = - {\vec {E} \left( {\nabla \hat
{\varepsilon }} \right)} \mathord{\left/ {\vphantom {{\vec {E}
\cdot \left( {\nabla \hat {\varepsilon }} \right)} {\hat
{\varepsilon }\,}}} \right. \kern-\nulldelimiterspace} {\hat
{\varepsilon }\,}$ from the equation $\nabla \vec {D} = 0$;
therefore, the second term is absent, both in (\ref{eq4}) written
in the time domain and in (\ref{eq6}), (\ref{eq9}), and
(\ref{eq14}), $\nabla \left( {\nabla \cdot \vec {E}}
\right) = 0$.

 According to \cite{[15]} and the
results of the previous section, one can show that the total field
$E_x \left( {y,z} \right) = E_x^{inc} \left( {y,z} \right) +
E_x^{scat} \left( {y,z} \right)$ of diffraction of the plane wave
\begin{equation}
\nonumber E_x^{inc} \left( {y,z} \right) = a^{inc} \exp \left\{
{\,i\,\left[ {\phi \;y - \Gamma  \left( {z - 2\,\pi \,\delta }
\right)} \right]\,} \right\}, \quad z > 2\,\pi \,\delta ,
\end{equation}
 by a weakly nonlinear dielectric layer
is the solution to the equation (see (\ref{eq4})):
\begin{equation}
\label{eq15}
 \nabla ^2 \cdot \vec {E} + \frac{\omega ^2}{c^2} \varepsilon ^{\left(
L \right)}\left( z \right) \cdot \vec {E} + \frac{4\pi \omega
^2}{c^2}
\vec {P}^{\left( {NL} \right)}
 \equiv \left[ {\nabla ^2 + \kappa ^2 \varepsilon \left( {z,\,\alpha
,\,\left| {E_x } \right|^2} \right)} \right]\,E_x \left( {y,z}
\right) = 0
\end{equation}
satisfying the following generalized boundary conditions:

\noindent continuity of $E_{tg}$ and $H_{tg}$ on the boundary of
the layer having the permittivity $\varepsilon \left( {z,\,\alpha
,\,\left| {E_x } \right|\,^2} \right)\mbox{ , }$

\noindent the spatial quasi-homogeneity condition \cite{[3]} with
respect to $y$
\begin{equation}
\label{eq16} E_x \left( {y,z} \right) = U\left( z \right) \exp
\left( {i\phi \;y} \right)\mbox{ ,}
\end{equation}

\noindent
and the radiation condition for the scattered field
\begin{equation}
\label{eq17}
E_x^{scat} \left( {y,z} \right) = \left\{ {{\begin{array}{*{20}c}
 {a^{scat}} \hfill \\
 {b^{scat}} \hfill \\
\end{array} }} \right\}\,\exp \left( {i {\kern 1pt} \left(
{\phi {\kern 1pt} \;y\pm \Gamma {\kern 1pt} \left( {z \mp 2\,\pi
\,\delta } \right)} \right)} \right){\kern 1pt} ,\;\quad
z{\begin{array}{*{20}c}
 > \hfill \\
 < \hfill \\
\end{array} }\pm 2\,\pi \,\delta .
\end{equation}


Here we use the following notations: $\{x,y,z,t\}$ are
dimensionless spatial-temporal coordinates introduced so that the
layer thickness is $4\pi \,\delta $; the time dependence is $\exp
\,( - \,i\,\omega \,t)$; $\omega =
 \kappa \,c$ is the dimensionless circular frequency; $\kappa = \omega
\mathord{\left/ {\vphantom {\omega c}} \right. \kern-\nulldelimiterspace} c
\equiv {2\pi } \mathord{\left/ {\vphantom {{2\pi } \lambda }} \right.
\kern-\nulldelimiterspace} \lambda $ is the dimensionless frequency
parameter such that ${h{\kern 1pt} } \mathord{\left/ {\vphantom {{h{\kern
1pt} } {\lambda \,}}} \right. \kern-\nulldelimiterspace} {\lambda \,} =
2{\kern 1pt} \kappa \,\delta $, where $\lambda $ is the free-space
wavelength; $c = (\varepsilon _0 \,\mu _0 )^{1 \mathord{\left/ {\vphantom {1
2}} \right. \kern-\nulldelimiterspace} 2}$ is the dimensionless quantity
equal to the speed of light in the medium containing the layer ($Imc = 0)$;
$\varepsilon _0 $ and $\mu _0 $ are the material parameters of the medium;
$E_{tg} $ and $H_{tg} $ are the tangential components of the total $\vec
{E}$ and $\,\vec {H}$ fields; $\nabla ^2 = {\partial \,^2} \mathord{\left/
{\vphantom {{\partial \,^2} {\partial {\kern 1pt} y\,^2}}} \right.
\kern-\nulldelimiterspace} {\partial {\kern 1pt} y\,^2} + {\partial \,^2}
\mathord{\left/ {\vphantom {{\partial \,^2} {\partial {\kern 1pt} \,z^2}}}
\right. \kern-\nulldelimiterspace} {\partial {\kern 1pt} \,z^2}$;

\begin{equation}
\nonumber
 \varepsilon \left( {z,\;\alpha ,\left| {E_x } \right|^2}
\right) \equiv \varepsilon \left( {z,\;\alpha ,\left| {U\left( z
\right)} \right|^2} \right) = \left\{ {{\begin{array}{*{20}c}
 {\quad \quad \;\;\;\quad \quad \quad \quad 1\;,\quad \left| z \right| >
2\pi \delta } \hfill \\
 {\varepsilon ^{\left( L \right)}\left( z \right) + \alpha  \left|
{U\left( z \right)} \right|^2,\quad \left| z \right| \le 2\pi \delta }
\hfill \\
\end{array} }} \right. \quad ,
\end{equation}

\noindent where $\varepsilon ^{\left( L \right)}\left( z \right)$
is piecewise continuously differentiable with respect to $z$;
$\alpha = 3\pi \chi _{xxxx}^{\left( 3 \right)} $; $\Gamma =
(\kappa \,^2 - \phi \,^2)\,^{1 \mathord{\left/ {\vphantom {1 2}}
\right. \kern-\nulldelimiterspace} 2}$ is the transverse
propagation constant (transverse wavenumber); $\phi \equiv \kappa
\cdot \sin \left( \varphi \right)$ is the longitudinal propagation
constant (longitudinal wavenumber); and $\varphi $ is the angle of
incidence of the plane wave, $\left| \varphi \right| < \pi
\mathord{\left/ {\vphantom {\pi 2}} \right.
\kern-\nulldelimiterspace} 2$. 
Quantities
${x}',{y}',{z}',{t}',{\omega }'$ are reconstructed from the
dimensionless values by the formulas $({x}',{y}',{z}') = (x,y,z)
\cdot h \mathord{\left/ {\vphantom {h {4\pi \,\delta }}} \right.
\kern-\nulldelimiterspace} {4\pi \,\delta }$, ${t}' = t \cdot h
\mathord{\left/ {\vphantom {h {4\pi \,\delta }}} \right.
\kern-\nulldelimiterspace} {4\pi \,\delta }$, and ${\omega }' =
\omega \, {4\pi \,\delta } \mathord{\left/ {\vphantom {{4\pi
\,\delta } {{\kern 1pt} h}}} \right. \kern-\nulldelimiterspace}
{{\kern 1pt} h}$.

We look for the solution to problem (\ref{eq15})--(\ref{eq17}) in
the form

\begin{equation}
\label{eq18}
\begin{array}{l}
 E_x \left( {y,z} \right) = U\left( z \right) \,\exp \left(
{i\,\;\phi \,\;y} \right) = \\
 = \left\{ {{\begin{array}{*{20}c}
 {a^{inc} \exp \left\{ {i \left[ {\phi y - \Gamma
\left( {z - 2\pi \delta } \right)} \right]} \right\} + a^{scat}
\exp \left\{ {i{\kern 1pt} \left[ {\phi y + \Gamma \left( {z -
2\pi \delta } \right)} \right]}
\right\},\;\;\;z > 2\pi\delta ,} \hfill \\
 {U^{scat}\left( z \right) \exp \left( {i{\kern 1pt} \phi
{\kern 1pt} y} \right),\,\quad \;\;\quad \,\quad \quad \quad \quad
\quad \quad \quad \quad \quad \;\,\quad \quad \quad \quad \quad
\quad \,\;\left| z
\right| \le 2\pi \delta ,} \hfill \\
 {b^{scat} \exp \left\{ {i{\kern 1pt} {\kern 1pt} \left[
{\phi y - \Gamma \left( {z + 2\pi\delta } \right)} \right]}
\right\},\,\;\;\quad \;\quad \,\quad \quad \quad \quad \quad \quad
\quad
\quad \quad \quad \quad z < - 2\pi \delta ,} \hfill \\
\end{array} }} \right. \\
 \end{array}
\end{equation}
assuming the continuity on the permittivity break lines $z = 2\pi
\delta $ and $z = - 2\pi \delta $, so that $U\left( { - 2\pi
\delta } \right) = b\,^{scat}$ and $U\left( {2\pi \delta } \right)
= a\,^{inc} + a\,^{scat}$.
%
%
\section{Integral equation of the nonlinear problem}
\subsection{Reducing to an integral equation }

We solve problem (\ref{eq15})--(\ref{eq17}) in the whole space $Q
= \left\{ {q = \left\{ {y,z} \right\}:\,\;\; - \infty < y,\,z <
\infty} \right\}$ by reducing it to a one-dimensional IE along the
layer height $z \in \left[ { - 2\pi \delta ,\,2\pi \delta }
\right]$ with respect to the scattered field component $U\left( z
\right) \equiv U^{scat}\left( z \right)$ introduced in
(\ref{eq18}). To this end, make use of canonical Green's function
$G_0 $ of problem (\ref{eq15})--(\ref{eq17}) (for $\varepsilon =
1)$ defined in the strip $Q_{\left\{ {Y,\infty } \right\}} =
\left\{ {q = \left\{ {y,z} \right\}:\; - Y \le y \le Y\,,\;\;\vert
z\vert < \infty ;\;Y
> 0} \right\} \subset Q$ by the expression \cite{[15]}, \cite{[23]}, \cite{[24]}
\begin{equation}
\label{eq19}
\begin{array}{l}
 G_0 \left( {q\,,q_0 } \right) = \frac{i}{4\,Y}\;{\exp \left\{ {i
\left[ {\phi \left( {y - y_0 } \right) + \Gamma \left| {z - z_0 }
\right|} \right]\,} \right\}\,} \mathord{\left/ {\vphantom {{\exp
\left\{ {i  \left[ {\phi \left( {y - y_0 } \right) + \Gamma \left|
{z - z_0 } \right|} \right]\,} \right\}\,} {\Gamma }}}
\right. \kern-\nulldelimiterspace} {\,\Gamma } \equiv \\
 \equiv \exp \left( {\pm i\phi y} \right)  \frac{i\pi
}{4Y}\int\limits_{ - \infty }^\infty {H_0^{\left( 1 \right)} \left[ {\kappa
\sqrt {\left( {\tilde {y} - y_0 } \right)^2 + \left( {z - z_0 } \right)^2} }
\right]\exp \left( { \mp i\phi \tilde {y}} \right)d\tilde {y}} \, .\\
 \end{array}
\end{equation}

The nonlinear IE with respect to $U\left( z \right)$ introduced in
(\ref{eq18}) is obtained using a classical approach set forth in
\cite{[25]}. Denote by $V\left( q \right) \equiv E_x \left( {q
\equiv \left\{ {y,\,z} \right\}} \right) = U\left( z \right)
\,\exp \left( {i\phi y} \right)$ the total diffraction field (see
(\ref{eq18}) ), where $U\left( z \right)$ is the solution of
problem (\ref{eq15})--(\ref{eq17}), and write equation
(\ref{eq15}) as

\begin{equation}
\label{eq20}
\left( {\nabla ^2 + \kappa ^2} \right)\,\,V\left( q \right) = \left[ {\,1 -
\varepsilon \left( {q,\,\alpha ,\,\left| {V\left( q \right)} \right|^2}
\right)} \right]^{\kern 1pt} \kappa ^2\;V\left( q \right).
\end{equation}

The field of the incident plane wave $V_0 \left( q \right) \equiv
V^{inc}\left( q \right) = a^{inc}  \exp \left\{ {i\left[ {\phi y -
\Gamma \left( {z - 2\pi\delta } \right)} \right]\,} \right\}\,$
satisfies in the whole space $Q$ the homogeneous Helmholtz
equation

\begin{equation}
\label{eq21}
\left( {\nabla ^2 + \kappa ^2} \right)\,\,V_0 \left( q \right) = 0.
\end{equation}

At $z > 2\pi \delta  $, $  V_0 \left( q \right)$ is the incident
field of the incoming plane wave irradiating the layer, while at
$z < 2\pi \delta  \quad V_0 \left( q \right)$ is the outgoing
plane wave that satisfies the radiation condition at infinity
(because in the representation for $V_0 \left( q \right)$ the
transverse wavenumber $\Gamma > 0)$.

Subtracting from (\ref{eq20}) equation (\ref{eq21}) we obtain

\begin{equation}
\label{eq22}
\left( {\nabla ^2 + \kappa ^2} \right)\,\,\left[ {V\left( q \right) - V_0
\left( q \right)} \right] = \left[ {\,1 - \varepsilon \left( {q,\,\alpha
,\,\left| {V\left( q \right)} \right|^2} \right)} \right]^{\kern 1pt} \kappa
^2\;V\left( q \right).
\end{equation}

Here $V\left( q \right) - V_0 \left( q \right)$ satisfies the
radiation condition (\ref{eq17}) in the whole space. In fact, at
$z > 2\pi \delta $ the difference $V\left( q \right) - V_0 \left(
q \right) = V^{scat}\left( q \right)$ is the reflected field and
at $z \to - \infty $ both $V\left( q \right)$ and $V_0 \left( q
\right)$ satisfies the radiation condition.

Using (\ref{eq22}) and the equation for the canonical Green's function $G_0 $

\begin{equation}
\label{eq23} \left( {\nabla ^2 + \kappa ^2} \right)\,G_0 \left(
{q,\,q_0 } \right) = - \delta \left( {q,\,q_0 } \right)
\end{equation}
(where $\delta \left( {q,\,q_0 } \right)$ is the Dirac
delta-function) it is easy to show that

\begin{equation}
\label{eq24} \left( {V - V_0 } \right) \nabla ^2G_0 - G_0 \nabla
^2\,\left( {V - V_0 } \right) = - \left( {V - V_0 } \right) \delta
\left( {q,\,q_0 } \right) - G_0  \left[ {1 - \varepsilon \left(
{q,\,\alpha ,\left| V \right|^2} \right)} \right]^{\kern 1pt}
\kappa ^2 V.
\end{equation}

Let $Q_{\left\{ {Y,\,Z} \right\}} = \left\{ {q = \left\{ {y,z}
\right\}:\,\; - Y \le y \le Y,\; - Z \le z \le Z;\;Y > 0,\,Z >
2\pi \delta } \right\}$ denote a rectangular domain in space $Q$.
Divide this domain into rectangles such that in each of them the
permittivity $\varepsilon \left( {q,\,\alpha ,\,\left| V
\right|^2} \right)$ is continuously differentiable with respect to
$y$ and $z$. On common parts of the boundaries of rectangles
$V\left( q \right)\,$ and ${\partial V\left( q \right)}
\mathord{\left/ {\vphantom {{\partial V\left( q \right)} {\partial
{\kern 1pt} n}}} \right. \kern-\nulldelimiterspace} {\partial
{\kern 1pt} n}$ (where $n$ denotes the outer normal) are
continuous due to the continuity of the tangential components
$E_{tg} \;$ and $H_{tg} $. Therefore, in the whole domain
$Q_{\left\{ {Y,\,Z} \right\}} $, the sought-for twice continuously
differentiable function $V\left( q \right)\,$ preserves this
property up to the boundary $\partial Q_{\left\{ {Y,\,Z} \right\}}
$; i.e., $V\left( q \right)\,\, \in C^2\left( {Q_{\left\{ {Y,\,Z}
\right\}} } \right) \cap C^1\left( {\overline {Q_{\left\{ {Y,\,Z}
\right\}} } } \right)$ (here $\overline {Q_{\left\{ {Y,\,Z}
\right\}} } = Q_{\left\{ {Y,\,Z} \right\}} \cup \partial
Q_{\left\{ {Y,\,Z} \right\}} )$.

Applying in $Q_{\left\{ {Y,\,Z} \right\}} $ Green's formula

\begin{equation}
\nonumber
 \int\!\!\!\int\limits_{Q_{\left\{ {Y,Z} \right\}} }
{\left[ {\left( {V - V_0 } \right) \nabla ^2G_0 - G_0  \nabla
^2\left( {V - V_0 } \right)} \right]\,dq_0 } =
\int\limits_{\partial Q_{\left\{ {Y,Z} \right\}} } {\left[ {\left(
{V - V_0 } \right) \frac{\partial G_0 }{\partial n} - G_0
\frac{\partial \left( {V - V_0 } \right)}{\partial n}} \right]dq_0
} ,
\end{equation}

\noindent
and taking into account (\ref{eq24}), we obtain

\begin{equation}
\label{eq25}
\begin{array}{l}
 V\left( q \right) = - \kappa ^2\;\int\!\!\!\int\limits_{Q_{\left\{ {Y,Z}
\right\}} } {G_0 \left( {q,q_0 } \right) \left[ {\,1 - \varepsilon
\left( {q_0 ,\,\alpha ,\,\left| {V\left( {q_0 } \right)}
\right|^2} \right)}
\right]^{\kern 1pt} V\left( {q_0 } \right)dq_0 } + V_0 \left( q \right) - \\
 \quad \quad \quad - \int\limits_{\partial Q_{\left\{ {Y,Z} \right\}} }
{\left\{ {\left[ {V\left( {q_0 } \right) - V_0 \left( {q_0 }
\right)} \right] \frac{\partial G_0 \left( {q,q_0 }
\right)}{\partial n} - G_0 \left( {q,q_0 } \right) \frac{\partial
\left[ {V\left( {q_0 } \right)
- V_0 \left( {q_0 } \right)} \right]}{\partial n}} \right\}dq_0 } \,. \\
 \end{array}
\end{equation}
When parameter $Z \to \infty $, the integrals in the lower,
$\left[ {\left( { - Z, - Y} \right),\,\left( { - Z,Y} \right)}
\right]$, and upper $\left[ {\left( {Z,Y} \right),\,\left( {Z, -
Y} \right)} \right]$ parts of the boundary $\partial Q_{\left\{
{Y,\,Z} \right\}} $ that enter curvilinear integral (\ref{eq25})
tend to zero. This statement follows from asymptotic properties of
Green's function (\ref{eq10}) and the fact that $V^{scat} = V -
V_0 $ satisfies radiation condition (\ref{eq17}). The integrals
along $\left[ {\left( { - Z,Y} \right),\,\left( {Z,Y} \right)}
\right]$ and $\left[ {\left( {Z, - Y} \right),\,\left( { - Z, - Y}
\right)} \right]$ cancel each other. Therefore, setting in
(\ref{eq25}) $Z \to \infty $ and deleting the curvilinear integral
along the boundary $\partial Q_{\left\{ {Y,\,Z \to \infty }
\right\}} $, we obtain an integral representation for the total
field of diffraction in the band $Q_{\left\{ {Y,\infty }
\right\}}$

\begin{equation}
\nonumber
 V\left( q \right) = - \kappa
^2\;\int\!\!\!\int\limits_{Q_{\left\{ {Y,\infty } \right\}} } {G_0
\left( {q,q_0 } \right) \left[ {\,1 - \varepsilon \left( {q_0
,\,\alpha ,\,\left| {V\left( {q_0 } \right)} \right|^2} \right)\;}
\right]^{\kern 1pt} V\left( {q_0 } \right)dq_0 } + V_0 \left( q
\right)\,,\,\quad q \in Q_{\left\{ {Y,\infty } \right\}} .
\end{equation}
The integrand in the remaining double integral is a finite
function with respect to $z$, i.e.
 $$
 1 - \varepsilon(
q,\alpha ,|V(q)|^2)\equiv 0,\quad \left| z \right| > 2\pi\delta ,
$$ 
because the permittivity of the medium enveloping the
layer is assumed to be equal unity, so that one can limit the
integration to the domain occupied by the dielectric

\begin{equation}
\nonumber V\left( q \right) = - \kappa ^2
\int\!\!\!\int\limits_{Q_{\left\{ {Y,Z = 2\pi \delta } \right\}} }
{G_0 \left( {q,q_0 } \right) \left[ {1 - \varepsilon \left( {q_0
,\alpha ,\left| {V\left( {q_0 } \right)} \right|^2} \right)}
\right]^{\kern 1pt} V\left( {q_0 } \right)dq_0 } + V_0 \left( q
\right),\quad q \in Q_{\left\{ {Y,\infty } \right\}} .
\end{equation}

Performing a transfer to the limit $Y \to \infty $ (which can be
justified by the facts that, according to (\ref{eq16}) and
(\ref{eq19}), the integrand is asymptotically equivalent to
$O\left( {Y^{ - 1}} \right)$ and parameter $Y$may be chosen
arbitrarily) we obtain an integral representation for the total
field of diffraction in the whole space $Q$

\begin{equation}
\label{eq26} V\left( q \right) = - \kappa
^2\;\int\!\!\!\int\limits_{Q_\delta } {G_0 \left( {q,q_0 } \right)
 \left[ {1 - \varepsilon \left( {q_0 ,\,\alpha ,\,\left| {V\left(
{q_0 } \right)} \right|^2} \right)} \right]^{\kern 1pt} V\left(
{q_0 } \right)dq_0 } + V_0 \left( q \right)\,,\quad q \in Q\,.
\end{equation}
Here $Q_\delta \equiv Q_{\left\{ {\infty ,Z = 2\pi \delta }
\right\}} = \left\{ {q = \left\{ {y,z} \right\}:\;\; - \infty < y
< + \infty \,,\;\;\vert z\vert \le 2\pi \,\delta } \right\}$ is
the band occupied by the nonlinear dielectric layer.

We can also obtain (\ref{eq26}) using an iteration scheme based on
the approach developed in \cite{[15]}, \cite{[24]}. Let us give a
short description of this method. In space $Q$ a function sequence
$V_n (y,z)$ is constructed such that every function of this
sequence, beginning from $n = 1$, satisfies conditions
(\ref{eq16}) and (\ref{eq17}), and the limiting function $V = E_x
\left( {y,z} \right) = \mathop {\lim }\limits_{n \to \infty }
\,V_n $ is a solution to (\ref{eq15})--(\ref{eq17}); namely,

\begin{equation}
\label{eq27}
\begin{array}{l}
 \left( {\nabla ^2 + \kappa ^2} \right)\,\,V_0 = 0\,,\quad \left( {\nabla ^2
+ \kappa ^2} \right)\,\,V_1 = \left[ {\,1 - \varepsilon \left( {z,\,\alpha
,\,\left| {V_0 } \right|^2} \right)} \right]^{\kern 1pt} \kappa ^2\;V_0 \, +
V_0 ,\dots \, ,\quad \\
 \left( {\nabla ^2 + \kappa ^2} \right)\,\,V_{n + 1} = \left[ {\,1 -
\varepsilon \left( {z,\,\alpha ,\,\left| {V_n } \right|^2} \right)}
\right]^{\kern 1pt} \,\kappa ^2\;V_n \, + V_0 ,...\quad \quad \quad \quad
\quad \quad \;\; \\
 \end{array}
\end{equation}

Equations (\ref{eq27}) are formally equivalent to the following

\begin{equation}
\label{eq28}
\begin{array}{l}
 V_0 \left( q \right) \equiv V^{ins}\left( q \right)\,,\\
 V_1
\left( q \right) = - \int\!\!\!\int\limits_{Q_\delta } {G_0 \left(
{q,\,q_0 } \right) \left[ {1 - \varepsilon \left( {q_0 ,\alpha
,\left| {V_0 \left( {q_0 } \right)} \right|^2} \right)}
\right]^{\kern 1pt} \kappa ^2V_0 \left( {q_0 } \right)d{\kern 1pt}
q_0 } + V_0 \left( q
\right)\, ,\,\dots  \, ,\\
 V_{n + 1} \left( q \right) = - \int\!\!\!\int\limits_{Q_\delta } {G_0
\left( {q,\,q_0 } \right) \left[ {\,1 - \varepsilon \left( {q_0
,\,\alpha ,\,\left| {V_n \left( {q_0 } \right)} \right|^2}
\right)} \right]^{\kern 1pt} \kappa ^2V_n \left( {q_0 }
\right)\,d{\kern 1pt} q_0 } +
V_0 \left( q \right)\,\dots  \, , \,\,\,\, q \in Q. \\
 \end{array}
\end{equation}

Performing in (\ref{eq28}) a transfer to the limit $n \to \infty $ we obtain the
integral representation (\ref{eq26}) for the total field of diffraction in $Q$.

For $q \in Q_\delta $, representation (\ref{eq26}) is transformed
to a nonlinear IE with respect to the sought for scattered field
$V\left( q \right) \equiv V^{scat}\left( q \right)\,$, $q \in
Q_\delta $, see (\ref{eq18}). Substituting into equation
(\ref{eq29}) formula (\ref{eq19}) for canonical Green's function
and the expression for the permittivity $\varepsilon \left( {q_0
,\,\alpha ,\,\;\left| {V\left( {q_0 \;} \right)} \right|\;^2}
\right)\;$ we obtain an equation
\begin{gather}
U\left( z \right)e^{i{\kern 1pt} \phi {\kern 1pt} y} =\notag \\
- \mathop {\lim }\limits_{Y \to \infty } \left\{ {\frac{i\,\kappa
\,^2}{4Y\,\Gamma }\int\limits_{ - 2\pi \delta }^{2\pi \delta }
{\,\int\limits_{ - Y}^Y {e^{i{\kern 1pt} \phi {\kern 1pt}
y}e^{i\,\,\Gamma \left| {z - z_0 } \right|}} \left[ {\,1 -
\,\left( {\varepsilon ^{\left( L \right)}\,\left( {z_0 } \right) +
\alpha \,\left| {U\left( {z_0 } \right)} \right|^2} \right)}
\right]^{\kern 1pt} U\left( {z_0 } \right)dy_0 dz_0 } }
\right\} \\
+ \,\,U^{inc}\left( z \right)e^{i{\kern 1pt} \phi {\kern 1pt} y}
\notag
\end{gather}
 with respect to $U\left( z \right)
\equiv U^{scat}\left( z \right)\,$, $\quad \left| z \right| \le
2\pi \delta \;$, that enters the expression for the field $V\left(
q \right) \equiv E_x \left( {q \equiv \left\{ {y,\,z} \right\}}
\right) = U\left( z \right) \exp \left( {i\phi y} \right)$
quasi-homogeneous along the layer. Integrating in domain $Q_\delta
$ with respect to $y_0 $, we obtain a nonlinear IE of the second
kind with respect to an unknown function $U\left( z \right) \in
L_2 \left( {\,\left[ { - 2\pi \delta ,\;2\pi \delta } \right]\,}
\right)$:
\begin{eqnarray}
\label{eq29}
U\left( z \right) &+& \frac{i\,\kappa \,^2}{2\,\Gamma
}\int\limits_{ - 2\pi \delta }^{2\pi \delta } {\exp \left(
{i\Gamma  \left| {z - z_0 } \right|} \right)}\left[ {1 - \left(
{\varepsilon ^{\left( L \right)}\left( {z_0 } \right) + \alpha
\left| {U\left( {z_0 } \right)} \right|^2} \right)} \right]{\kern
1pt}U\left( {z_0 } \right)\,dz_0  \\
\nonumber &=& U^{inc}\left( z \right),\,\quad \left| z \right| \le
2\pi \delta ,
\end{eqnarray}
where $U^{inc}\left( z \right) = a^{inc}\exp \,\left[ { - i{\kern
1pt} \Gamma \left( {z - 2\pi\delta } \right)} \right]$.

The existence and uniqueness of the solution to IE (\ref{eq29})
for the linear problem with $\alpha = 0$ are proved in
\cite{[15]}, \cite{[26]}. In the general case a nonlinear IE of
type (\ref{eq29}) may or may not have (the unique) solution. Its
solvability is governed by the properties of the kernel and the
right-hand side (incident field $U^{inc}\left( z \right))$ and
value of the nonlinearity parameter.

Note that from the method of obtaining IE (\ref{eq29}) it follows
that the solution to this IE may be used for the integral
representation in $E_x \left( {y,z} \right) = U\left( z \right)
\,\exp \left( {i\phi y} \right)$ of the sought for solution to
problem (\ref{eq15})--(\ref{eq17}) for the points with the
coordinates $\left| z \right| > 2\pi \delta $ outside the
nonlinear layer. Indeed, finding the solution to IE (\ref{eq29})
and substituting it under the integral sign in (\ref{eq29}), we
obtain an explicit expression for $U\left( z \right)$ in the
domain $\left| z \right| > 2\pi \delta $.


The equivalence of IE (\ref{eq29}) to problem
(\ref{eq15})--(\ref{eq18}) is proved in Appendix.

\subsection{Sufficient condition of the existence of solution to nonlinear IE}

Assume that the permittivity function $\varepsilon ^{\left( L
\right)}\,\left( {z_0 } \right)$ is positive, bounded, and continuous in the
interval $\gamma = \left[ { - 2\pi \delta ,\,2\pi \delta } \right]$, so that
$\varepsilon ^{\left( L \right)}\,\left( {z_0 } \right) \in C\left( \gamma
\right)$, where $C\left( \gamma \right)$ denotes the space of continuous
functions in the closed interval $\gamma $ with the norm $\left\| f \right\|
= \left\| f \right\|_{C\left( \gamma \right)} = \mathop {\max }\limits_{z
\in \gamma } \left| {f\left( z \right)} \right|$. Assume also that

\begin{equation}
\label{eq30}
1 < \varepsilon ^{(L)}\left( z \right) \le E\,,\quad z \in \gamma \,,\quad E
> 1.
\end{equation}

Write IE (\ref{eq29}) in the operator form

\begin{equation}
\label{eq31}
U + AU - \alpha \,{\kern 1pt} \,F\left( U \right) = f\;,
\end{equation}

\noindent
where

\begin{equation}
\nonumber
 f\left( z \right) = U^{inc}\left( z \right) = a\exp
\left[ { - i\,\kappa \,\,\cos \left( \varphi \right) \left( {z -
2\pi \delta } \right)} \right]\,,\quad a = a^{inc} > 0\,,
\end{equation}

\begin{equation}
\label{eq32} AU = \int\limits_{ - 2\pi \delta }^{2\pi \delta }
{k\left( {z - z_0 } \right)  \left[ {1 - \,\varepsilon ^{\left( L
\right)}\,\left( {z_0 } \right)} \right]\,{\kern 1pt} U\left( {z_0
} \right)\,dz_0 }
\end{equation}

\noindent is a linear integral operator with the continuous kernel

\begin{equation}
\nonumber k\left( t \right) = s_0 \exp \left[ {2\,i\,\kappa
\,\,\cos \left( \varphi \right) \left| {\,t\,} \right|}
\right]\,,\quad s_0 = \frac{i\,\kappa \,}{2\,\cos \left( \varphi
\right)}\,\quad \left( { - \frac{\pi }{2\,} < \varphi < \frac{\pi
}{2\,}} \right),
\end{equation}
and
\begin{equation}
\nonumber F\left( U \right) = \int\limits_{ - 2\pi \delta }^{2\pi
\delta } {k\left( {z - z_0 } \right) \left| {U\left( {z_0 }
\right)} \right|^2U\left( {z_0 } \right)\,\,dz_0 }
\end{equation}
is a cubic-nonlinear integral operator.

A linear integral operator $B:\;C\left( \gamma \right) \to C\left( \gamma
\right)$ defined by
\begin{equation}
\nonumber Bu = \int\limits_{ - 2\pi \delta }^{2\pi \delta }
{k\left( {z - z_0 } \right) u\left( {z_0 } \right)\,\,dz_0 } \;
\end{equation}
is bounded, so that
\begin{equation}
\label{eq33}
\left\| B \right\| \le \int\limits_{ - 2\pi \delta }^{2\pi \delta } {\mathop
{\max }\limits_{z \in \gamma } \left| {k\left( {z - z_0 } \right)}
\right|dz_0 } = 4\pi \delta \frac{\kappa }{2\cos \left( \varphi \right)} =
\frac{2\pi \delta \kappa }{\cos \left( \varphi \right)} = q_0 .
\end{equation}

Integral operator (\ref{eq33}) is bounded and continuous in
$C\left( \gamma \right)$ and its norm can be estimated as
\begin{equation}
\label{eq34} \left\| A \right\| \le \mathop {\max }\limits_{z \in
\gamma } \left[ {\int\limits_{ - 2\pi \delta }^{2\pi \delta }
{\left| {k\left( {z - z_0 } \right)} \right| \left| {1 -
\varepsilon ^{\left( L \right)}\left( {z_0 } \right)} \right|dz_0
} } \right] = \left( {E - 1} \right)\frac{2\pi \delta \kappa
}{\cos \left( \varphi \right)} = \left( {E - 1} \right)q_0 .
\end{equation}
The nonlinear operator $Q\left( U \right) = \left| U \right|^2
 U$ is bounded and continuous in $C\left( \gamma \right)$ and
$F\left( U \right)$ is therefore bounded and continuous in
$C\left( \gamma \right)$ as a superposition of $B$ and $Q$. Hence
the nonlinear operator $T\left( U \right) = - AU + \alpha \,{\kern
1pt} \,F\left( U \right) + f\;$is completely continuous on each
bounded subset $\Omega \in C\left( \gamma \right)$.

Set
$$\Omega = S_p = \left\{ {U \in C\left( \gamma
\right):\;\;\left\| U \right\| < p} \right\}$$ to be a ball in
$C\left( \gamma \right)$, assume that $U \in S_p $, and estimate
the $C\left( \gamma \right)$-norm of $T$:
\[
\begin{array}{l}
 \left\| {T\left( U \right)} \right\| \le \left\| A \right\|  \left\| U
\right\| + \alpha \,\left\| B \right\|\,{\kern 1pt}  \left\| U
\right\|^3\, + \left\| f \right\| \le \\
 \quad \quad \;\;\; \le p\,\left\| A \right\| + \alpha \,\left\| B
\right\|\,p^3 + a \le K\left( p \right)\,,\quad \quad U \in S_p \,, \\
 K\left( p \right) = \left( {E - 1} \right)q_0 p + \alpha \,q_0 p^3 + a\,.\,
\\
 \end{array}
\]

Write equation (\ref{eq31}) as $U = T\left( U \right)$. One can
apply to the operator $T\left( U \right)\,:\;S_p \to S_p $ the
Banach fixed-point theorem \cite{[27]} if $K\left( p \right) \le
p\;$. In order to determine the corresponding range of values of
parameter $p > 0$, solve the inequality $K\left( p \right) \le
p\;$, which yields a cubic inequality with respect to $p$
\begin{equation}
\label{eq35}
P_K \left( p \right) \equiv D_0 p^3 - D_1 p + a \le 0,\quad \quad D_0 =
\alpha \,q_0 ,\quad \quad D_1 = 1 - \left( {E - 1} \right)\,q_0 .
\end{equation}

The necessary condition for (\ref{eq35}) to have a positive solution is $D_1 > 0$,
which yields $\left( {E - 1} \right)\,q_0 < 1$, or, according to (\ref{eq33}) and
assumption (\ref{eq30}) concerning the properties of the permittivity function
$\varepsilon ^{\left( L \right)}\left( z \right)$,
\begin{equation}
\label{eq36} 2\pi \delta \kappa < \cos \left( \varphi \right)
\left\{ {\mathop {\max }\limits_{z \in \gamma } \left[
{\varepsilon ^{\left( L \right)}\left( z \right)} \right] - 1}
\right\}^{ - 1}\;\quad \quad \left( { - \frac{\pi }{2\,} < \varphi
< \frac{\pi }{2\,}} \right).
\end{equation}

Subject to the condition (\ref{eq36}), cubic polynomial $P_K
\left( p \right)$ in (\ref{eq35}) has two positive zeros $p_1 $
and $p_2 $ if $0 < a < \mathop {\max }\limits_{p \ge 0} ( - D_0
p^3 + D_1 p)$ (with $D_1 > 0)$, which holds if the local minimum
of $P_K (p)$ is negative, $P_K \left( {p_{ext} } \right) < 0$, at
the point $p_{ext} = \sqrt {\frac{D_1 }{3D_0 }} > 0$ where ${P}'_K
\left( {p_{ext} } \right) = 0$. The corresponding condition for
$a$ can be written as $a < \frac{2}{3}D_1 \sqrt {\frac{D_1 }{3D_0
}} $, or
\begin{equation}
\label{eq37}
a\sqrt \alpha < \frac{2}{3\sqrt 3 }\frac{\left[ {1 - \left( {E - 1}
\right)\,q_0 } \right]^{3 / 2}}{\sqrt {q_0 } }\,,\quad \quad q_0 =
\frac{2\pi \delta \kappa }{\cos \left( \varphi \right)}.
\end{equation}

Condition (\ref{eq37}) holds for arbitrary set of the problem parameters $a,\kappa
,\varphi ,\delta ,$ and $E$ satisfying (\ref{eq36}) if the nonlinearity parameter
$\alpha $ is sufficiently small because $\alpha $ enters only the left-hand
side of inequality (\ref{eq37}). The inequality $K\left( p \right) \le p\;$ holds
for$p \in \left( {p_1 ,p_2 } \right)$, $p_1 > p_2 > 0$; for example, at

\begin{equation}
\nonumber p = p_{ext} = \sqrt {\frac{1 - \left( {E - 1}
\right)\,q_0 }{3\,\alpha \,q_0 }} .
\end{equation}


\begin{theorem}
Assume that

\noindent (i)\quad  the permittivity function $\varepsilon
^{\left( L \right)}\left( z \right)$ is positive, bounded, and
continuous in the closed interval $\gamma = \left[ { - 2\pi \delta
,\,2\pi \delta } \right]$ and $E = \mathop {\max }\limits_{z \in
\gamma } \left[ {\varepsilon ^{\left( L \right)}\left( z \right)}
\right] > 1$;

\noindent (ii)\quad  parameters $a$, $\kappa$, $\delta $ and the
nonlinearity parameter $\alpha $ are positive, $\left| {\,\varphi
\,} \right| < \pi / 2$, and all the parameters $a$, $\kappa $,
$\delta $, $\varphi $, $\alpha $, and $E$ satisfy (\ref{eq30}),
(\ref{eq36}) and (\ref{eq37}), namely

\begin{equation}
\nonumber \quad \quad \quad \quad \quad E > 1,\quad \left( {E - 1}
\right)q_0 < 1\,,\quad \alpha < \alpha _0 =
\frac{4}{27}\frac{1}{a^2}\frac{\left[ {\,1 - \left( {E - 1}
\right)q_0 } \right]^{{\kern 1pt} 3}}{q_0 }\,,\quad q_0 =
\frac{2\pi \delta \kappa }{\cos \left( \varphi \right)}\,.
\end{equation}

Then the operator $T(U) = - AU + \alpha \,{\kern 1pt} \,F(U) + f$,
$T\left( U \right)\,:\;S_p \to S_p $ defined by (\ref{eq31}) and (\ref{eq32}) is a
contraction in the space $C\left( \gamma \right)$ if

\begin{equation}
\label{eq38} t_0 = q_0 \left( {E - 1 + 3\,\alpha \,p^2} \right) <
1\,,\quad p \in \left( {p_1 ,p_2 } \right)
\end{equation}
where $p_1 $ and $p_2$ are positive zeros of the polynomial
$P_K(p)$ defined by (\ref{eq35}) and $\alpha $ is sufficiently
small, satisfying
\begin{equation}
\label{eq39}
0 < \alpha < \min \left\{ {\alpha _0 ,\alpha _1 } \right\}\,,\quad \alpha _1
= \frac{1}{3}\left\{ {q_0 \left[ {\,1 - \left( {E - 1} \right)q_0 }
\right]\,} \right\}^{ - 1}.
\end{equation}

\end{theorem}

\textbf{Proof}.
Use definition (\ref{eq31}), estimates (\ref{eq33}) and
(\ref{eq34}), and inequality $\left| {\,\left| {z_1 } \right| -
\left| {z_2 } \right|\,} \right| \le \left| {z_1 - z_2 } \right|$
(where $z_1 $, $z_2 $ are complex numbers), assume that $U,V \in
S_p $, and estimate the $C\left( \gamma \right)$-norm of the
difference $T\left( U \right) - T\left( V \right)$:

\begin{equation}
\nonumber
\begin{array}{l}
 \left\| {T\left( U \right) - T\left( V \right)} \right\| \le \left\| {AV -
AU} \right\| + \alpha \,\,\left\| {F\left( U \right) - F\left( V \right)}
\right\| \le \\
 \quad \quad \quad \quad \quad \;\; \le \left( {E - 1} \right)q_0 \left\| {U
- V} \right\| + \alpha \,q_0 \left\| {U\left| U \right|^2 - V\left| V
\right|^2} \right\| < q_0 \left( {E - 1 + 3\alpha \,p^2} \right)\,\left\| {U
- V} \right\|\,. \\
 \end{array}
\end{equation}

Thus, inequality (\ref{eq38}) provides that operator $T\left( U
\right)\,:\,S_p \to S_p $ is a contraction if $\alpha $ is
sufficiently small; namely, satisfies (\ref{eq37}) and
(\ref{eq39}). Note that (\ref{eq39}) follows from (\ref{eq37}),
the condition $3\alpha {\kern 1pt} p_{ext}^2 < 1$, and inequality
$0 < p_1 < p_{ext} $, where $p_{ext} \in \left( {p_1 ,p_2 }
\right)$ is the point of a negative local minimum of the cubic
polynomial $P_K \left( p \right)$ (\ref{eq35}) satisfying ${P}'_K
(p_{ext} ) = 0$ and $\mathop {\min }\limits_{t \ge 0} \left[ {P_K
\left( t \right)} \right] = P_K \left( {p_{ext} } \right) < 0$ and
$p_1 $, $p_2 $ are positive zeros of $P_K \left( p \right)$.

Summarizing the results verified above we conclude that the
following statement is valid.


\begin{theorem}
Assume that the permittivity function $\varepsilon ^{\left( L
\right)}\left( z \right)$, parameters $a$, $\kappa $,$\delta $,
the nonlinearity parameter $\alpha $, and quantity $E$ satisfy
conditions (i) and (ii) from Theorem 1 and conditions (\ref{eq38})
and (\ref{eq39}). Then the operator $T\left( U \right) = - AU +
\alpha \,{\kern 1pt} F\left( U \right) + f$ defined by
(\ref{eq31}) and (\ref{eq32}) is a contraction in the space
$C\left( \gamma \right)$ and IE (\ref{eq29}) has the unique
solution $U^\ast \left( z \right)$ continuous in the closed
interval $\left[ { - 2\pi \delta ,\,2\pi \delta } \right]$.
$U^\ast \left( z \right)$ is a limit with respect to the $C\left(
\gamma \right)$-norm of the function sequence $U_n \left( z
\right)$ (the fixed point of operator $T\left( U \right))$
determined according to

\begin{equation}
\label{eq40} U_{n + 1} = T\left( {U_n } \right),\quad n =
0,\,1,\,2,\,\dots \,,\quad U_0 \in S_p = \left\{ {U \in C\left(
\gamma \right):\;\left\| U \right\| < p} \right\}\,,\quad p \in
(p_1 ,p_2 ),
\end{equation}

\noindent where $p_1 $ and $p_2 $ are positive zeros of the
polynomial $P_K \left( p \right)$ defined by (\ref{eq35}).
\end{theorem}

The rate of convergence of the fixed-point iterations (\ref{eq40})
can be estimated using the quantity $t_0 < 1$ defined in
(\ref{eq38}):

\begin{equation}
\nonumber \left\| {U_n - U^\ast } \right\| = \left\| {T\left(
{U_{n - 1} } \right) - T\left( {U^\ast } \right)} \right\| < t_0
\left\| {U_{n - 1} - U^\ast } \right\| < \ldots < t_0 ^{n -
1}\left\| {T\left( {U_0 } \right) - U^\ast } \right\|,\quad n =
2,3,\dots \, .
\end{equation}


\subsection{Complex-valued permittivity function (diffraction
by a lossy nonlinear layer)}

The method and results can be extended to the case when the permittivity
$\varepsilon ^{\left( L \right)}\left( z \right)$ is an arbitrary
complex-valued function of the real argument $z$ continuous and bounded on
the line. To this end denote
\begin{equation}
\label{eq41}
\varepsilon ^{\left( L \right)}\left( z \right) - 1 = g\left( z \right) =
\varepsilon _1 \left( z \right)\exp \left[ {i\varepsilon _2 \left( z
\right)} \right] = g_1 \left( z \right) + i\,g_2 \left( z \right),
\end{equation}
where, according to physical assumptions of the model, the real and
imaginary parts of the permittivity function, $g_1 \left( z \right)$ and
$g_2 \left( z \right)$, are positive, continuous, and bounded on the line
satisfying $g_1 \left( z \right) \ge 1$ and $g_1 \left( z \right) > > g_2
\left( z \right)$, so that the modulus $\varepsilon _1 \left( z \right)$ and
argument $\varepsilon _2 \left( z \right)$ of the permittivity are also
positive functions continuous and bounded on the line with $0 \le
\varepsilon _2 \left( z \right) < \pi / 2$.

Make use of (\ref{eq41}) and represent integral operator (\ref{eq32}) as

\begin{equation}
\label{eq42}
A_1 U = \int\limits_{ - 2\pi \delta }^{2\pi \delta } {k_1 \left( {z,z_0 }
\right)\varepsilon _1 \left( {z_0 } \right){\kern 1pt} \,U\left( {z_0 }
\right)\,\,dz_0 } \,,\quad k_1 \left( {z,z_0 } \right) = - s_0 \exp \left\{
{\,i\,\left[ {2\kappa \cos \left( \varphi \right)\,\left| {z - z_0 }
\right|\,\varepsilon _2 \left( {z_0 } \right)} \right]} \right\}.
\end{equation}

Assuming, similar to (\ref{eq30}) and taking into account (\ref{eq41}) and the conditions
for the permittivity function, that

\begin{equation}
\label{eq43}
0 < \varepsilon _1 \left( z \right) \le E_1 \,,\quad z \in \gamma ,
\end{equation}
(that is, $0 < \left| {\varepsilon ^{\left( L \right)}\left( z
\right)} \right| \le E_1 $, $z \in \gamma )$ we can estimate, as
in (\ref{eq32}), the norm of the integral operator (\ref{eq42}),
which is bounded and continuous in $C\left( \gamma \right)$, as

\begin{equation}
\label{eq44} \left\| {A_1 } \right\| \le \mathop {\max }\limits_{z
\in \gamma } \left[ {\int\limits_{ - 2\pi \delta }^{2\pi \delta }
{\left| {k_1 \left( {z,z_0 } \right)} \right| \left| {\varepsilon
_1 \left( {z_0 } \right)} \right|dz_0 } } \right] = E_1 \frac{2\pi
\delta \kappa }{\cos \left( \varphi \right)} = E_1 q_0 .
\end{equation}

(\ref{eq44}) yields an estimate for the norm of the nonlinear operator $T_1 \left( U
\right) = - A_1 \,U + \alpha \,\,F\left( U \right) + f\;$

\begin{equation}
\label{eq45}
\left\| {T_1 \left( U \right)} \right\| \le K_1 \left( p \right)\,,\quad U
\in S_p \,,\quad K_1 \left( p \right) = E_1 q_0 p + \alpha \,q_0 p^3 + a.
\end{equation}

Thus one can easily check that the following statements are valid
which are extensions of Theorems 1 and 2 to the case of a
complex-valued permittivity function.


\begin{theorem}
Assume that

\noindent (i)\quad  the permittivity $\varepsilon ^{\left( L
\right)}\left( z \right)$ is a complex-valued function given by
(\ref{eq41}), where $g_1 \left( z \right) > 0$ and
%
%
%
are continuous and bounded on the line so that the modulus $\varepsilon _1
(z)$ and argument $\varepsilon _2 (z)$ of the function $\varepsilon ^{\left(
L \right)}\left( z \right) - 1$ are also nonnegative functions continuous
and bounded on the line with $0 \le \varepsilon _2 \left( z \right) < \pi /
2$ and $\varepsilon _1 \left( z \right) \ge 1$, and $E_1 = \mathop {\max
}\limits_{z \in \gamma } \left[ {\varepsilon _1 \left( z \right)} \right] >
0$ in the closed interval $\gamma = \left[ { - 2\pi \delta ,\,2\pi \delta }
\right]$;

\noindent (ii)\quad  parameters $a$, $\kappa $, $\delta $ and the
nonlinearity parameter $\alpha $ are positive, $\left| \varphi
\right| < \pi / 2,$ and all the parameters $a$, $\kappa $, $\delta
$, $\varphi $, $\alpha $, and $E_1 $ satisfy the conditions
similar to (\ref{eq30}), (\ref{eq36}) and (\ref{eq37}),  namely,

\begin{equation}
\label{eq46}
0 < E_1 q_0 < 1\,,\quad \alpha < \alpha _0 ^{(\ref{eq1})} =
\frac{4}{27}\frac{1}{a^2}\frac{\left( {1 - E_1 q_0 } \right)^3}{q_0
}\,,\quad q_0 = \frac{2\pi \delta \kappa }{\cos \left( \varphi \right)} >
0.
\end{equation}

Then the operator $T_1 \left( U \right) = - A_1 U + \alpha
\,\,F\left( U \right) + f$, $T_1 \left( U \right)\,\;:\;S_p \to S_p $
defined using (\ref{eq42}) is a contraction in the space$C\left( \gamma \right)$ if

\begin{equation}
\label{eq47}
t_1 = q_0 \left( {E_1 + 3\alpha {\kern 1pt} p^2} \right) < 1\,,\quad p \in
\left( {p_1 ^{\left( 1 \right)},p_2 ^{\left( 1 \right)}} \right),
\end{equation}

\noindent where $p_1 ^{\left( 1 \right)}$and $p_2 ^{\left( 1
\right)}$ are positive zeros of the polynomial

\begin{equation}
\label{eq48}
P_K ^{\left( 1 \right)}\left( p \right) \equiv D_0 p^3 - D_1 ^{\left( 1
\right)}p + a\,,\quad D_0 = \alpha \,q_0 ,\quad D_1 ^{(\ref{eq1})} = 1 - E_1 \,q_0
\end{equation}

\noindent and $\alpha $ is sufficiently small, satisfying

\begin{equation}
\label{eq49}
0 < \alpha < \min {\begin{array}{*{20}c}
 {\{\alpha _0 ^{(\ref{eq1})},\alpha _1 ^{(\ref{eq1})}\},} \hfill & {\alpha _1 ^{(\ref{eq1})}} \hfill
\\
\end{array} } = \frac{1}{3}\left[ {q_0 \left( {1 - E_1 q_0 } \right)}
\right]^{\, - 1}.
\end{equation}

\end{theorem}

\begin{theorem}

Assume that the permittivity function $\varepsilon ^{\left( L
\right)}\left( z \right)$ specified by (\ref{eq41}), parameters
$a$, $\kappa $, $\delta $, the nonlinearity parameter $\alpha $,
and quantity $E_1 = \mathop {\max }\limits_{z \in \gamma } \left[
{\varepsilon _1 \left( z \right)} \right] > 0$, $\gamma = \left[ {
- 2\pi \delta ,\,2\pi \delta } \right]$ ($\varepsilon _1 \left( z
\right) = \left| {\varepsilon ^{(L)}\left( z \right) - 1}
\right|)$ satisfy conditions (i) and (ii) from Theorem 3 and
conditions (\ref{eq46}), (\ref{eq47}), and (\ref{eq49}). Then the
operator $T_1 \left( U \right) = - A_1 U + \alpha \,\,F\left( U
\right) + f$ defined using (\ref{eq42}) is a contraction in the
space $C\left( \gamma \right)$ and IE (\ref{eq29}) has the unique
solution $U^\ast \left( z \right)$ continuous in the closed
interval $\left[ { - 2\pi \delta ,\,2\pi \delta } \right]$.
$U^\ast \left( z \right)$ is a limit with respect to the $C\left(
\gamma \right)$-norm of the function sequence $U_n \left( z
\right)$ (the fixed point of operator $T_1 \left( U \right))$
determined according to

\begin{equation}
\label{eq50}
U_{n + 1} = T_1 \left( {U_n } \right)\,\;,\quad n = 0,\;1,\;2,\;...,\quad
U_0 \in S_p = \left\{ {\,U \in C\left( \gamma \right):\;\left\| U \right\| <
p} \right\}\,,\quad p \in \left( {p_1 ^{\left( 1 \right)},\;p_2 ^{\left( 1
\right)}} \right),
\end{equation}

\noindent
where $p_1 ^{\left( 1 \right)}$and $p_2 ^{\left( 1 \right)}$are positive
zeros of the polynomial $P_K ^{\left( 1 \right)}\left( p \right)$ defined by
(\ref{eq48}).
\end{theorem}

The rate of convergence of the fixed-point iterations (\ref{eq50})
can be estimated using the quantity $t_1 < 1$ defined in
(\ref{eq47}):

\begin{equation}
\nonumber
\begin{array}{l}
 \left\| {\,U_n - U^\ast } \right\| < t_1 ^{n - 1}\left\| {\,T_1 \left( {U_0
} \right) - U^\ast } \right\|\,,\quad n = 2,\,3,\,\dots \, , \\
 U_0 \in S_{p^\ast } = \left\{ {\,U \in C\left( \gamma \right):\;\left\| U
\right\| < p^\ast } \right\}\,,\quad p^\ast \in \left( {p_1 ^{\left( 1
\right)},p_2 ^{\left( 1 \right)}} \right)\,. \\
 \end{array}
\end{equation}

The existence of the unique solution to IE (\ref{eq29}) subject to
the sufficient conditions specified in formulations of Theorems 2
and 4 (corresponding to the cases of, respectively, real- and
complex-valued permittivity function of the nonlinear layer) and
the equivalence of IE (\ref{eq29}) to problem
(\ref{eq15})--(\ref{eq18}) proved in Appendix enables us to prove
the following statement which constitutes the main result of this
study.

\begin{theorem}
Assume that the permittivity function $\varepsilon ^{\left( L
\right)}\left( z \right)$ is (a) real-valued and quantity $E$and
parameters $a = a^{inc}$, $\kappa $, $\delta $, $\varphi $, and
$\alpha $ satisfy conditions (i) and (ii) from Theorem 1,
(\ref{eq38}), and (\ref{eq39}); or (b) complex-valued (given by
(\ref{eq41})) and quantity $E_1 = \mathop {\max }\limits_{z \in
\gamma } \left| {\varepsilon ^{(L)}\left( z \right) - 1} \right| >
0$ and parameters $a = a^{inc}$, $\kappa $, $\delta $, $\varphi $,
and $\alpha $ satisfy conditions (i) and (ii) from Theorem 3 and
conditions (\ref{eq46}), (\ref{eq47}), and (\ref{eq49}). Then
problem (\ref{eq15})--(\ref{eq18}) has the unique solution $U^\ast
\left( z \right)$ continuous in the closed interval $\left[ { -
2\pi \delta ,\,2\pi \delta } \right]$ which can be determined as a
limit with respect to the $C\left( \gamma \right)$-norm of the
function sequence $U_n \left( z \right)$ determined, respectively,
according to (a) (\ref{eq40}) or (b) (\ref{eq50}).
\end{theorem}

\subsection{First iterations as trigonometric polynomials}
In view of the fact that at $\varepsilon ^{\left( L \right)}\left(
z \right) = 1$ and $\alpha = 0$ nonlinear IE (\ref{eq20}) has a
formal solution
\begin{equation}
\label{eq51}
U\left( z \right) = U^{inc}\left( z \right) = \tilde {a}\exp \left( { -
i\,bz} \right)\,,\quad \tilde {a} = \exp \left( {ibd} \right)\,,\quad d =
2\pi \delta \,,\quad b = \kappa \,\cos \left( \varphi \right),
\end{equation}
it is reasonable to choose the zero iteration in (\ref{eq50}) $U_0 \left( z \right)
= \tilde {a}\exp \left( { - i\,bz} \right)$ in the form (\ref{eq51}). The linear
integral operators (\ref{eq32}) and (\ref{eq42}) can be represented as
\begin{equation}
\label{eq52}
\begin{array}{l}
 AU = A\left[ \eta \right]\,U = \int\limits_{ - d}^d {k\left( {z - z_0 }
\right)\,\eta \left( {z_0 } \right)U\left( {z_0 } \right)\,\,dz_0 } \,, \\
 k\left( t \right) = s_0 \exp \left( {2i\,b\left| {\,t\,} \right|}
\right)\,,\quad s_0 = \frac{i\,\kappa \,}{2\,\cos \left( \varphi
\right)}\,,\quad \left( { - \frac{\pi }{2\,} < \varphi < \frac{\pi }{2\,}}
\right)\,.\; \\
 \end{array}
\end{equation}

Obviously, they are linear with respect to the (continuous complex-valued)
weight function $\eta \left( {z_0 } \right)$:

\begin{equation}
\nonumber A\left[ {h_1 \eta _1 + h_2 \eta _2 } \right]\,U = h_1
A\left[ {\eta _1 } \right]\,U + h_2 A\left[ {\eta _2 }
\right]\,U\,\;,\quad h_1 ,\,h_2 = const.
\end{equation}

\begin{lemma}

\begin{equation}
\nonumber
A\left[ {\eta ^{(0)}} \right]\,U_0 = H_1 \exp \left( { -
i\,2bz} \right) + H_2 \exp \left( {i\,2bz} \right) + H_3 \exp
\left[ {i\,z\left( {q - b} \right)} \right],
\end{equation}

\noindent
where

\begin{equation}
\nonumber
H_j = H_0 \tilde {H}_j \,,\quad \left( {j = 1,\;2,\;3}
\right)\,,\quad H_0 = - \frac{iT\,\tilde {a}s_0 }{\left( {b + q}
\right)\left( {3b - q} \right)},
\end{equation}

\begin{equation}
\label{eq53}
\tilde {H}_1 = \left( {3b - q} \right)\exp \left[ {i\,d\left( {b + q}
\right)} \right]\,,\quad \tilde {H}_2 = \left( {b + q} \right)\exp \left[
{i\,d\left( {3b - q} \right)} \right]\,,\quad \tilde {H}_3 = - 4b,
\end{equation}

$q \ne - b$, $q \ne 3b,$ and

\begin{equation}
\label{eq54}
U_0 \left( z \right) = \tilde {a}\exp \left( { - i\,bz} \right)\,,\quad \eta
^{\left( 0 \right)}\left( z \right) = T\exp \left( {i\,qz} \right)\,,\quad
\tilde {a},\;T = const.
\end{equation}

At $\quad q = b$,

\begin{equation}
\nonumber
A\left[ {\eta ^{\left( 0 \right)}} \right]\,U_0 = -
\frac{iT\,\tilde {a}s_0 }{b}\left[ {\exp \left( {i\,2bd}
\right)\cos \left( {i2bz} \right) - 1} \right].
\end{equation}

\end{lemma}

\textbf{Proof}. Proof of Lemma 1 reduces to tedious algebra and
integration.

We see that it is possible to determine explicitly the image $A\left[ {\eta
^{\left( 0 \right)}} \right]^{\kern 1pt} \,U_0 $ of a simple trigonometric
polynomial $U_0 = \tilde {a}\exp \left( { - i\,bz} \right)$and to show that
this image is also a trigonometric polynomial:

\begin{equation}
\nonumber
A\left[ {\eta ^{\left( 0 \right)}} \right]\,\left[
{\tilde {a}\exp \left( { - i\,bz} \right)} \right] =
\sum\limits_{j = 1}^3 {H_j \exp \left( {i\,c_j z} \right)}
\,,\quad c_1 = - 2b\,,\quad c_2 = 2b\,,\quad c_3 = q - b.
\end{equation}

The linearity of $A\left[ \eta \right]\,U$ with respect to the
weight function $\eta \left( z \right)$ and $U$ yields

\begin{lemma}

Let

\begin{equation}
\nonumber
\eta ^{\left( 0 \right)}\left( z \right) =
\sum\limits_{j = 1}^{N_\eta } {r_j \exp \left( {i\,q_j z} \right)}
\,,\quad N_\eta \ge 1.
\end{equation}

Then the image $A\left[ {\eta ^{\left( 0 \right)}} \right]\,U$ of
a trigonometric polynomial $U\left( z \right) = \sum\limits_{j =
1}^{N_U } {h_j \exp \left( {i\,b_j z} \right)} $ is also a
trigonometric polynomial:

\begin{equation}
\nonumber
A\left[ {\eta ^{\left( 0 \right)}} \right]\;U\left( z
\right) = \sum\limits_{j = 1}^{N_A } {P_j \exp \left( {i\,c_j z}
\right)} ,
\end{equation}

\noindent where the coefficients $P_j $ and the number of terms
$N_A$ can be determined explicitly.

\end{lemma}

Similar statements are valid for the nonlinear operator $F\left( U \right)$ defined
in (\ref{eq31}); namely, the image $F(U_0 )$ of a trigonometric polynomial $U_0 =
\tilde {a}\exp \left( { - i\,bz} \right)$ is also a trigonometric polynomial
that can be determined explicitly:


\begin{lemma}

\begin{equation}
\nonumber
F\left( {U_0 } \right) = f_1 \exp \left( { - i\,2bz}
\right) + f_2 \exp \left( {i\,2bz} \right) + f_3 \exp \left( { -
i\,bz} \right),
\end{equation}

\noindent where

\begin{equation}
\nonumber
\begin{array}{l}
 f_j = f_0 \tilde {f}_j \,,\quad \left( {j = 1,\;2,\;3} \right)\,,\quad f_0
= - \frac{i\,a^3\exp \left( {i{\kern 1pt} b{\kern 1pt} d} \right)s_0
}{3b}\,, \\
 \tilde {f}_1 = 3\exp \left( {i{\kern 1pt} db} \right)\,,\quad \tilde {f}_2
= - \exp \left( {i{\kern 1pt} d{\kern 1pt} b} \right)\,,\quad \tilde {f}_3 =
- 2\,, \\
 U_0 \left( z \right) = \tilde {a}\exp \left( { - i\,bz} \right)\,,\quad
\tilde {a} = a\exp \left( {i\,bd} \right)\,,\quad a = const\,. \\
 \end{array}
\end{equation}

\end{lemma}

Lemmas 1--3 enable one to evaluate explicitly the first iteration

\begin{equation}
\nonumber
\begin{array}{l}
 U_1 = T\left( {U_0 } \right) = - AU_0 + \alpha \,\,F\left( {U_0 } \right) +
U_0 = \; \\
 \quad \; = - \sum\limits_{j = 1}^3 {H_j \exp \left( {i\,c_j z} \right)} +
\alpha \sum\limits_{j = 1}^2 {f_j \exp \left( {i\,c_j z} \right)} + \alpha
{\kern 1pt} f_3 \exp \left( { - i\,bz} \right) = \\
 \quad \; = \sum\limits_{j = 1}^4 {S_j \exp \left( {i\,c_j z} \right)}
\,,\quad S_j = - H_j + \alpha {\kern 1pt} f_j \quad \left( {j = 1,\,2}
\right)\,,\quad S_3 = \alpha {\kern 1pt} f_3 \,,\quad S_4 = - H_3 \,, \\
 \end{array}
\end{equation}

\begin{equation}
\nonumber c_1 = - 2b\,,\quad c_2 = 2b\,,\quad c_3 = - b\,,\quad
c_4 = q - b\,.
\end{equation}

We conclude that according to (\ref{eq32}) and (\ref{eq52}) if the permittivity function
$\varepsilon ^{\left( L \right)}\left( z \right)$ is a trigonometric
polynomial then the first iteration $U_1 $ specified by (\ref{eq50}) is also a
trigonometric polynomial whose coefficients can be determined explicitly.

\subsection{Sufficient condition of the existence of solution to nonlinear IE:
reducing to a functional equation system}
Here we present the
proof of an alternative sufficient condition for the existence of
a solution to nonlinear IE (\ref{eq29}) which is similar to the
solvability conditions of the type (\ref{eq38}). The approach
developed in this section enables one to create a rather efficient
method of the numerical solution of the IE. To this end, reduce
(\ref{eq29}) to a nonlinear functional equation system,
considering the system of two IEs in the domain $\quad \left| z
\right| \le 2\pi \delta \;$:

\begin{equation}
\label{eq55}
\begin{array}{l}
 U_{n + 1} \left( z \right) + \frac{i\kappa^2}{2\Gamma }\int\limits_{
- 2\pi \delta }^{2\pi \delta } {\exp \left( {i\Gamma \left| {z -
z_0 } \right|} \right) \left[ {1 - \left( {\varepsilon^{\left( L
\right)}\,\left( {z_0 } \right) + \alpha \,\left| {U_n \left( {z_0
} \right)} \right|^2} \right)} \right]\,{\kern 1pt} U_n \left(
{z_0 }
\right)\,dz_0 } = U^{inc}\left( z \right), \\
 \Psi _n \left( z \right) + \frac{i\kappa^2}{2\Gamma }\int\limits_{ -
2\pi \delta }^{2\pi \delta } {\exp \left( {i\Gamma \left| {z - z_0
} \right|} \right) \left[ {1 - \left( {\varepsilon^{\left( L
\right)}\left( {z_0 } \right) + \alpha \left| {U_n \left( {z_0 }
\right)} \right|^2} \right)} \right]\,{\kern 1pt} \Psi _n \left(
{z_0 } \right)\,\,dz_0 } = U^{inc}\left( z \right). \\
 \end{array}
\end{equation}

The first equation of system (\ref{eq55}) is an iteration scheme
of solution to nonlinear equation (\ref{eq29}) (cf. (\ref{eq28})).
The second is a linear IE with respect to $\Psi _n \left( z
\right)$ for the given $U_n \left( {z_0 } \right)\,$. If $\Psi _n
\left( z \right)\;$ is not an eigenfunction of the problem of
diffraction by the layer with the permittivity $\varepsilon \left(
{z,\,\alpha ,\,\left| {U_n \left( z \right)} \right|^2} \right)
\equiv \varepsilon ^{\left( L \right)}\,\left( z \right) + \alpha
\,\left| {U_n \left( z \right)} \right|^2$, then the second
equation is uniquely solvable [15, 26] and its solution can be
represented in the form

\begin{equation}
\label{eq56} \Psi _n \left( z \right) = \Psi \left( {z,\,\alpha
,\,\left| {U_n \left( z \right)} \right|^2} \right) U^{inc}\left(
z \right).
\end{equation}

\noindent
where $\Psi \left( {z,\,\alpha ,\,\left| {U_n \left( z \right)} \right|^2}
\right)$ is the solution to the linear IE at $U^{inc}\left( z \right) = 1$
such that $\left| {\Psi \left( {z,\,\alpha ,\,\left| {U_n \left( z \right)}
\right|^2} \right)\,} \right| \le 1$.

The analysis of the convergence criterion for the sequence $U_n \left( z
\right)\,$, $\Psi _n \left( z \right)\;$ specified by system (\ref{eq55}) enables
one to obtain a sufficient condition for the existence of solution to
nonlinear IE (\ref{eq29}).

Kernels of IEs (\ref{eq55}) are identical, which makes it possible to calculate and
estimate the $L_2 $-norm of the difference between $U_n \left( z
\right)\,$and $\Psi _n \left( z \right)\;$

\begin{equation}
\label{eq57}
\begin{array}{l}
 \rho \left[ {U_{n + 1} \left( z \right),\,\Psi _n \left( z \right)} \right]
= \left[ {\int\limits_{ - 2\pi \delta }^{2\pi \delta } {\left| {U_{n + 1}
\left( z \right) - \Psi _n \left( z \right)} \right|^2dz\,} } \right]^{1
\mathord{\left/ {\vphantom {1 2}} \right. \kern-\nulldelimiterspace} 2} = \\
 = \left| {\frac{i\,\kappa \,^2}{2\,\Gamma }} \right| \left\{
{\int\limits_{ - 2\pi \delta }^{2\pi \delta } {\left|
{\int\limits_{ - 2\pi \delta }^{2\pi \delta } {\exp \left(
{i\Gamma \left| {z - z_0 } \right|} \right) } \left[ {1 - \left(
{\varepsilon ^{\left( L \right)}\left( {z_0 } \right) + \alpha
\left| {U_n \left( {z_0 } \right)} \right|^2} \right)}
\right]{\kern 1pt}  \left[ {U_n \left( {z_0 } \right) - \Psi _n
\left( {z_0 } \right)} \right]dz_0 } \right|^2dz} } \right\}^{1
\mathord{\left/ {\vphantom {1 2}} \right.
\kern-\nulldelimiterspace} 2} = \\
 = \frac{\kappa \,^2}{2\,\Gamma } \left\{ {\int\limits_{ - 2\pi \delta
}^{2\pi \delta } {\left| {\int\limits_{ - 2\pi \delta }^{2\pi
\delta }  \left[ {1 - \left( {\varepsilon ^{\left( L
\right)}\left( {z_0 } \right) + \alpha\left| {U_n \left( {z_0 }
\right)} \right|^2} \right)} \right]{\kern 1pt} \left[ {U_n \left(
{z_0 } \right) - \Psi _n \left( {z_0 } \right)} \right]dz_0 }
\right|^2dz} } \right\}^{1 \mathord{\left/ {\vphantom {1 2}}
\right. \kern-\nulldelimiterspace} 2} \le
\\
 \le \frac{\kappa^2}{2\,\Gamma } \left\{ {\int\limits_{ - 2\pi
\delta }^{2\pi \delta } \, \int\limits_{ - 2\pi \delta }^{2\pi
\delta } \left| {1 - \left( {\varepsilon ^{\left( L \right)}\left(
{z_0 } \right) + \alpha \left| {U_n \left( {z_0 } \right)}
\right|^2} \right)} \right|^2dz_0\, dz} \right\}^{1
\mathord{\left/ {\vphantom {1 2}} \right.
\kern-\nulldelimiterspace} 2}  \left\{ {\int\limits_{ - 2\pi
\delta }^{2\pi \delta }  \left| {U_n \left( {z_0 } \right) - \Psi
_n \left( {z_0 } \right)} \right|^2dz_0 } \right\}^{1
\mathord{\left/
{\vphantom {1 2}} \right. \kern-\nulldelimiterspace} 2} \le \\
 \le \frac{\kappa \,^2}{2\,\Gamma } 4\pi \delta \mathop {\max
}\limits_{\left| z \right| \le 2\pi \delta } \left| {1 - \left(
{\varepsilon ^{\left( L \right)}\left( z \right) + \alpha\left|
{U_n \left( z \right)} \right|^2} \right)} \right| \rho \left[
{U_n
\left( z \right),\,\Psi _n \left( z \right)} \right] \le \\
 \, \le \frac{\kappa \,^2}{2\,\Gamma }  4\pi \delta \mathop {\max
}\limits_{\left| z \right| \le 2\pi \delta } \left[ {\left| {1 -
\varepsilon ^{\left( L \right)}\left( z \right)} \right| + \left|
\alpha \right| \left| {U_n \left( z \right)} \right|^2} \right]
\rho
\left[ {U_n \left( z \right),\,\Psi _n \left( z \right)} \right] \le \\
 \le \frac{\kappa \,^2}{2\,\Gamma }  4\pi \delta  \mathop {\max
}\limits_{\left| z \right| \le 2\pi \delta } \left[ {\left| {1 -
\varepsilon ^{\left( L \right)}\left( z \right)} \right| + \left|
\alpha \right| \left| {U^{inc}\left( z \right)} \right|^2} \right]
\rho
\left[ {U_n \left( z \right),\,\Psi _n \left( z \right)} \right]\,. \\
 \end{array}
\end{equation}
The last inequality in (\ref{eq57}) is obtained taking into account the condition
$\mathop {\max }\limits_{\left| z \right| \le 2\pi \delta } \left| {U_n
\left( z \right)} \right| \le \mathop {\max }\limits_{\left| z \right| \le
2\pi \delta } \left| {U^{inc}\left( z \right)} \right|$ which holds for all
$n = 0,\,1,\,2,\,\ldots $ and directly follows from the inequality $\mathop
{\max }\limits_{\left| z \right| \le 2\pi \delta } \left| {U\left( z
\right)} \right| \le \mathop {\max }\limits_{\left| z \right| \le 2\pi
\delta } \left| {U^{inc}\left( z \right)} \right|$ due to (\ref{eq18}). We see that,
according to (\ref{eq57}), in the case under study of weakly nonlinear approximation
(\ref{eq5}) when

\begin{equation}
\label{eq58} \mathop {\max }\limits_{\left| z \right| \le 2\pi
\delta } \left[ {\,\left| \alpha \right|  \left| {U\left( z
\right)} \right|^2} \right] \le \mathop {\max }\limits_{\left| z
\right| \le 2\pi \delta } \left[ {\left| \alpha \right| \left|
{U^{inc}\left( z \right)} \right|^2} \right] < \mathop {\max
}\limits_{\left| z \right| \le 2\pi \delta } \left| {\varepsilon
^{\left( L \right)}\,\left( z \right)} \right|,
\end{equation}

\noindent the iterations defined by the first equation of
(\ref{eq55}) converge to the unique solution determined by
(\ref{eq55}) if the term in the last inequality of (\ref{eq57})
multiplying the norm satisfies the condition

\begin{equation}
\nonumber \frac{\kappa \,^2}{2\,\Gamma }  4\pi \delta \mathop
{\max }\limits_{\left| z \right| \le 2\pi \delta } \left[ {\left|
{\,1 - \varepsilon ^{\left( L \right)}\,\left( z \right)} \right|
+ \left| \alpha \right|  \left| {U^{inc}\left( z \right)}
\right|^2} \right] < 1.
\end{equation}

\noindent In view of the expression for the transverse wavenumber
$\Gamma = \left\{ {\kappa ^2 - \left[ {\kappa \sin \left( \varphi
\right)} \right]^{{\kern 1pt} 2}} \right\}^{{\kern 1pt} 1
\mathord{\left/ {\vphantom {1 2}} \right.
\kern-\nulldelimiterspace} 2}$ rewrite the last inequality as

\begin{equation}
\label{eq59} \kappa \,  2\pi \delta \mathop {\max }\limits_{\left|
z \right| \le 2\pi \delta } \left[ {\left| {\,1 - \varepsilon
^{\left( L \right)}\,\left( z \right)} \right| + \left| \alpha
\right| \left| {U^{inc}\left( z \right)} \right|^2} \right] < \cos
\left( \varphi \right).
\end{equation}

Note that according to (\ref{eq58}), condition (\ref{eq59}) can be written in the form (\ref{eq38})
with $p = \frac{a}{\sqrt 3 }$ as

\begin{equation}
\label{eq60} q_0 \left( {E - 1 + \alpha \,a^2} \right) < 1\,.
\end{equation}

Observe also that (\ref{eq59}) yields the sufficient condition
(\ref{eq35}) for the existence of solution to nonlinear IE
(\ref{eq29}) obtained in the previous section because the latter
reads $\left( {E - 1} \right)  q_0 < 1$.

We have proved the following statement which constitutes a sufficient
condition for the existence of solution to nonlinear IE (\ref{eq29}).

\begin{theorem}
Assume that the weakly nonlinear approximation
(\ref{eq58}) holds. Then nonlinear IE (\ref{eq29}) has the unique
continuous solution if condition (\ref{eq59}) holds. This solution
can be obtained using both the iterations defined by the first
equation of (\ref{eq55}) and the equivalent iteration scheme
according to the second equation of (\ref{eq55}) if to consider
its solution $\Psi _n \left( z \right)\;$ as the $n + 1$
approximation (setting $\Psi _n \left( z \right)\;\, \equiv U_{n +
1} \left( z \right))$ to the sought for $U\left( z \right)\,$.
\end{theorem}
\section{Conclusion}
%
%
We have proved, subject to certain sufficient
conditions, the unique solvability of the problem of diffraction
of a plane wave by a transversely inhomogeneous isotropic
nonmagnetic linearly polarized dielectric layer filled with a
Kerr-type nonlinear medium. The diffraction problem has been
reduced to a cubic-nonlinear IE of the second kind. Based on the
use of the contraction principle, sufficient conditions of the IE
unique solvability have been obtained in the form of simple
inequalities. The method presented in this work can be generalized
so that it will enable one to obtain eigensolutions and
soliton-type solutions; eigenvalues, also as functions of the
problem parameters; and to develop the techniques to wider classes
of nonlinearities $B$ and operators $L$ of singular semilinear
BVPs $L\left( \lambda \right)\,u + \alpha B\left( {u{\kern 1pt}
;\lambda } \right) = f$associated with the problems of wave
scattering and propagation$. $

One the basis of these solution techniques and the IE obtained one
can perform numerical investigation of the resonance effects
caused by certain nonlinear properties of the object under study
irradiated by an intense electromagnetic field. In particular, one
can determine the critical limits of the excitation field
intensity that govern applicability of the developed mathematical
model. The proposed methods and results of computations can be
further applied to the analysis of various physical phenomena
including self-influence and interaction of waves; determination
of eigenfields, natural (resonance) frequencies of nonlinear
objects, and dispersion amplitude--phase characteristics of the
diffraction fields; description of evolution processes in the
vicinities of critical points; and to the design and modeling of
novel scattering, transmitting, and memory devices.

\bigskip
\noindent {\large{\bf{Appendix }}}
\bigskip

\noindent Let us prove that IE (\ref{eq29}) is equivalent to BVP
(\ref{eq15})--(\ref{eq18}); namely, if $U\left( z \right)$ is a
solution to IE (\ref{eq29}) then $E_x \left( {y,z} \right) =
U\left( z \right) \,\exp \left( {i\,\phi \,y} \right)$ is a
solution to (\ref{eq15})--(\ref{eq17}) subject to representation
(\ref{eq18}) and \textit{vice versa}. To this end, let us show
that IE (\ref{eq29}) and (\ref{eq15})--(\ref{eq18}) are reduced to
the determination of the solution to one and the same BVP and both
problems are equivalent to one and the same IE. Indeed, write IE
(\ref{eq29}) for the points $\left| z \right| \le 2\pi \delta $
inside the nonlinear layer in the form
\begin{equation}
\label{eq61}
\begin{array}{l}
 U\left( z \right) + \frac{i\,\kappa \,^2}{2\,\Gamma }[F_+ (z) + F_ -
(z)] = U^{inc}\left( z \right),\quad \quad \left| z \right| \le
2\pi
\delta \;, \\
 \,F_ + \,\,(z) = \int\limits_{ - 2\pi \delta }^z {\exp \left[ {i\,\Gamma
 \left( {z - z_0 } \right)} \right]} \left[ {\,1 - \,\left(
{\varepsilon ^{\left( L \right)}\,\left( {z_0 } \right) + \alpha
\,\left| {U\left( {z_0 } \right)} \right|^2} \right)}
\right]\,{\kern 1pt}
U\left( {z_0 } \right)\,\,dz_0 , \\
 F_ - \,\,(z) = \int\limits_z^{2\pi \delta } {\exp \left[ { - i\,\Gamma
 \left( {z - z_0 } \right)} \right]}  \left[ {\,1 - \,\left(
{\varepsilon ^{\left( L \right)}\,\left( {z_0 } \right) + \alpha
\,\left| {U\left( {z_0 } \right)} \right|^2} \right)}
\right]\,{\kern 1pt}
U\left( {z_0 } \right)\,\,dz_0 \,. \\
 \end{array}
\end{equation}
We have
\begin{equation}
\label{eq62}
\begin{array}{l}
 U\left( {2\pi \delta } \right) = U^{inc}\left( {2\pi \delta } \right) -
\frac{i\,\kappa \,^2}{2\,\Gamma }F_ + \,\,(2\pi \delta ) = a^{inc} +
a^{scat}, \\
 U\left( { - 2\pi \delta } \right) = U^{inc}\left( { - 2\pi \delta } \right)
- \frac{i\,\kappa \,^2}{2\,\Gamma }F_ - \,\,( - 2\pi \delta ) =
a^{inc}e^{4i\pi \delta } - \frac{i\,\kappa \,^2}{2\,\Gamma }F_ - \,\,( -
2\pi \delta ) = b^{scat}, \\
 \end{array}
\end{equation}
where
\begin{equation}
\label{eq63}
a^{scat} = - \frac{i\,\kappa \,^2}{2\,\Gamma }F_ + \,\,(2\pi \delta
),{\begin{array}{*{20}c}
 \hfill & \hfill & {b^{scat} = a^{inc}e^{4i\pi \delta } - \frac{i\,\kappa
\,^2}{2\,\Gamma }F_ - \,\,( - 2\pi \delta )} \hfill \\
\end{array} }
\end{equation}
denote the quantities (constants) expressed in terms of the solution to IE
(\ref{eq29}).

Differentiating two times with respect to $z$ the first equality
(\ref{eq61}) (or IE (\ref{eq29})) involving functions $F_\pm
\,(z)$ and using the continuity condition (\ref{eq16}) for the
tangential components of the total diffraction field on the
permittivity break lines $z = 2\pi \delta $ and $z = - 2\pi \delta
$ and representation (\ref{eq18}) we arrive at the problem in the
differential form equivalent to IE (\ref{eq61}) (or (\ref{eq29}))
\begin{eqnarray}
\label{eq63-1}
\begin{array}{l}
 \ell _\Gamma (U) + g(U) \equiv {U}''\left( z \right) + \left\{ {\Gamma ^2
- \kappa ^2\left[ {1 - \left( {\varepsilon ^{\left( L
\right)}\left( z \right) + \alpha \left| {U\left( z \right)}
\right|^2} \right)} \right]} \right\}^{\kern 1pt} U\left( z
\right) = 0,\quad  \left| z \right| \le
2\pi \delta , \\
 {\begin{array}{*{20}c}
 {{\begin{array}{*{20}c}
 {{\begin{array}{*{20}c}
 \hfill & \hfill \\
\end{array} }} \hfill & \hfill \\
\end{array} }} \hfill & \hfill & \hfill \\
\end{array} }U\left( {2\pi \delta } \right) = a^{scat} + a^{inc},\quad
U\left( { - 2\pi \delta } \right) = b^{scat},\quad \quad \;\;\; \\
 \end{array}
\end{eqnarray}
where the linear differential and nonlinear
operators
\begin{equation}
\nonumber
\ell _\Gamma (U) = {U}''\left( z \right) + \Gamma^2
U\left( z \right),{\begin{array}{*{20}c}
 \hfill & \hfill & {g(U) = - \kappa ^2\left[ {1 - \,\left( {\varepsilon
^{\left( L \right)}\,\left( z \right) + \alpha \,\left| {U\left( z \right)}
\right|^2} \right)} \right]^{\kern 1pt} U\left( z \right)} \hfill \\
\end{array} }.
\end{equation}

Note that $U(z)$ in (\ref{eq63-1}) satisfies the conditions
\begin{equation}
\label{eq64}
{U}'\left( {2\pi \delta } \right) = i\Gamma \left( {a^{scat} - a^{inc}}
\right),\quad \quad \;\;\;{U}'\left( { - 2\pi \delta } \right) = - i\Gamma
b^{scat},\quad
\end{equation}
Indeed, differentiating with respect to $z$ the first equality
(\ref{eq61}) and setting $z = 2\pi \delta $ and $z = - 2\pi \delta
$ we obtain
\begin{equation}
\nonumber
\begin{array}{l}
 {U}'\left( {2\pi \delta } \right) + i\,\,\Gamma \frac{i\,\kappa
\,^2}{2\,\Gamma }F_ + \,\left( {2\pi \delta } \right)\, = - i\,\,\Gamma
a^{inc}, \\
 {U}'\left( { - 2\pi \delta } \right) - i\,\,\Gamma \frac{i\,\kappa
\,^2}{2\,\Gamma }F_ - \,\left( { - 2\pi \delta } \right)\, = - i\,\,\Gamma
a^{inc}\exp \left( {4i\Gamma \pi \delta } \right), \\
 \end{array}
\end{equation}
which, together with (\ref{eq62}) and (\ref{eq63}), leads to
(\ref{eq64}). The same result is obtained if we note that on the
lines $z = 2\pi \delta $ and $z = - 2\pi \delta $, $U\,\,(z)$ and
its derivative coincide, according to (\ref{eq18}) and the
continuity condition, with the boundary values on these lines of
the respective functions $U_ + \left( z \right) = a^{inc} \exp
\left\{ {\, - i\Gamma \left( {z - 2\,\pi \,\delta } \right)\,}
\right\} + a^{scat} \exp \left\{ {\,i{\kern 1pt} \,\Gamma \left(
{z - 2\,\pi \,\delta } \right)\,} \right\}$, $U_ - \left( z
\right) = b^{scat}\exp \left\{ {\, - i{\kern 1pt} \,\Gamma \left(
{z + 2\,\pi \,\delta } \right)\,} \right\}$ and their derivatives.

Excluding in (\ref{eq63-1}) complex amplitudes $a^{scat}$ and
$b^{scat}$ we obtain a semilinear BVP of the Sturm--Liouville type
\begin{equation}
\label{eq65}
\begin{array}{l}
 {U}''\left( z \right) + \left\{ {\,\Gamma ^2 - \kappa ^2\left[ {1 -
\,\left( {\varepsilon ^{\left( L \right)}\,\left( z \right) + \alpha
\,\left| {U\left( z \right)} \right|^2} \right)} \right]\,} \right\}^{\kern
1pt} U\left( z \right) = 0,\quad \quad \left| z \right| \le 2\pi \delta \;,
\\
 i\Gamma U\left( {2\pi \delta } \right) - {U}'\left( {2\pi \delta } \right)
= 2i\Gamma a^{inc}, \\
 i\Gamma U\left( { - 2\pi \delta } \right) + {U}'\left( { - 2\pi \delta }
\right) = 0. \\
 \end{array}
\end{equation}

We can show independently that (\ref{eq65}) is equivalent to IE
(\ref{eq29}). Indeed, using Green's function
$$G \left( z,z_0
\right) = \frac{i}{2\,\Gamma }e^{i\Gamma \left| {z - z_0 }
\right|}$$ of the linear differential operator $\ell _\Gamma (U)$
we obtain IE (\ref{eq61}) (or (\ref{eq29})) by inverting $\ell
_\Gamma (U)$ in (\ref{eq63-1}) with the help of Green's function
$G \left( {z,z_0 } \right)$ (that is, by reducing (\ref{eq65}) to
an equivalent IE \cite{[28]}). The solution to IE (\ref{eq29})
satisfies the boundary condition of BVP (\ref{eq65}) at $z = \pm
2\pi \delta $ which is verified directly, as above, by
differentiating with respect to $z$ and setting $z = 2\pi \delta
$and $z = - 2\pi \delta $.

Assume now that $U\left( z \right)$ is a solution to IE
(\ref{eq29}) continuous in the closed interval $\left| z \right|
\le 2\pi \delta $ (we note that the unique solvability of IE
(\ref{eq29}) is proved in Section 4 subject to the sufficient
conditions formulated in Theorems 1--4). Then constants
$a\,^{scat}$and $b\,^{scat}$ are determined from (\ref{eq62}) and
(\ref{eq63}) so that $U\left( z \right)$ satisfies boundary
conditions in (\ref{eq63-1}) and (\ref{eq65}), and the solution to
(\ref{eq15})--(\ref{eq17}) is represented in the form (\ref{eq18})
with these constants (which also enter boundary conditions in
(\ref{eq65})).

The same BVP (\ref{eq65}) in the interval $\left| z \right| \le
2\pi \delta $ is obtained from the initial BVP
(\ref{eq15})--(\ref{eq17}) and representation (\ref{eq18}). This
follows directly if we substitute $E_x \left( {y,z} \right) =
U\left( z \right) \exp \left( {i\,\phi \,y} \right)$ into equation
(\ref{eq15}) taking into account the relationship $\Gamma ^2 =
\kappa ^2 - \phi ^2$ and the continuity of the tangential
components of the total diffraction field on the permittivity
break lines.

This statement completes the proof of the fact that IE
(\ref{eq29}) is equivalent to BVP (\ref{eq15})--(\ref{eq18}).

\bigskip


\label{lastpage}

\end{document}